\documentclass[a4paper,10pt]{amsart}
\usepackage[english]{babel} 
\usepackage[leqno]{amsmath}
\usepackage{amssymb,amsthm}
\usepackage{amscd}
\usepackage{enumerate}

\newtheorem{thm}{Theorem}[subsection]
\newtheorem{lemma}[thm]{Lemma}
\newtheorem{lemmadef}[thm]{Lemma - Definition}
\newtheorem{prop}[thm]{Proposition}

\newtheorem{cor}[thm]{Corollary}

\newtheorem{fact}[thm]{Fact}
\theoremstyle{remark}
\newtheorem{remark}[thm]{Remark}
\theoremstyle{definition}
\newtheorem{defi}[thm]{Definition}
\newtheorem{nota}[thm]{}

\newtheorem{example}[thm]{Example}

\newcommand{\la}{\longrightarrow}
\newcommand{\ha}{\hookrightarrow}

\newcommand{\ov}{\overline}

\newcommand{\Supp}{\operatorname{Supp}}

\newcommand{\im}{\operatorname{Im}}
\newcommand{\codim}{\operatorname{codim}}

\newcommand{\Pic}{\operatorname{Pic}}
\newcommand{\supp}{\operatorname{Supp}}
\newcommand{\mdeg}{\operatorname{{\underline{de}g}}}

\def\O{\mathcal O}

\def\X{\mathcal X}

\newcommand{\Z}{\mathbb{Z}}

\def\mc{\underline{c}}
\def\me{\underline{e}}
\def\md{\underline{d}}
\def\mt{\underline{t}}

\newcommand{\g}{\gamma}
\newcommand{\gS}{\gamma_S}
\newcommand{\dS}{\delta_S}
\newcommand{\dd}{\delta}
\newcommand{\nS}{\nu_S}
\newcommand{\eSd}{\epsilon^{\underline d}_S}
\newcommand{\nT}{\nu_T}
\newcommand{\gT}{\gamma_T}
\newcommand{\YT}{Y_T}
\newcommand{\dT}{\delta_T}

\newcommand{\rW}{\rho}

\newcommand{\pr}[1]{\mathbb{P}^{#1}}

\newcommand{\mgbar}{\ov{M}_g}

\def\aXd{\alpha^d_X}

\newcommand{\dfib}{\dot{X}}
\newcommand{\W}{{\operatorname{W}}}

\def\PXb{\overline{P^d_X}}

\def\PXg{P^{g-1}_X}
\def\PSd{P^{\md}_S}
\def\TSd{\Theta^{\md}_S}
\def\pXgb{\overline{P^{g-1}_X}}
\def\PXgb{\overline{P^{g-1}_X}}

\def\pfb{\overline{P^d_f}}

\def\Wmd{W_{\md}(X)}
\def\Wmdr{W_{\md}^r(X)}
\def\WMr{W_M^r(X)}

\def\BN{\rho(g,r,d)}

\def\Amd{A_{\md}(X)}
\def\Xmd{X^{\md}}
\def\dXmd{\dot{X}^{\md}}
\def\amd{\alpha^{\md}_X}
\def\WM{W_M(X)}

\def\FM{F_{M}(X)}

\def\Lv{(L_1,\ldots,L_{\gamma})}
\def\jv{j=1,\ldots,\delta}

\def\NE{N^{E_T}}

\def\WMS{W_M(X,S)}

\def\BXs{\Sigma^{ss} (X)}
\def\BX{\Sigma(X)}
\def\BY{\Sigma (Y)}
\def\BYs{\Sigma^{ss} (Y)}
\def\BYS{\Sigma (Y_S)}
\def\BYSs{\Sigma^{ss} (Y_S)}
\def\Sd{S(\underline{d})}
\def\Yd{Y_S(\underline{d})}
\def\nd{\nu_{\underline{d}}}
\def\dst{\underline{d}^{st}}
\def\sumg{\sum_1^{\gamma}g_i}

\newcommand{\picX}[1]{\Pic^{#1}X}

\newcommand{\YS}{Y_S}
\newcommand{\hX}{{\widehat{X}}}
\newcommand{\hXS}{{\widehat{X}}_S}

\newcommand{\hM}{{\widehat{M}}}

\newcommand{\sing}{X_{\text{sing}}}
\newcommand{\sep}{X_{\text{sep}}}

\begin{document}
\begin{center}
{\bf\large Geometry of the Theta Divisor of a compactified Jacobian}

\

\noindent {{Lucia  Caporaso}}\footnote{Dip.to di Matematica, Universit\`a Roma Tre,
L.go S.L.Murialdo 1,  00146 Roma-Italy, caporaso@mat.uniroma3.it}

\end{center}

%\bigskip \today

\

%\tableofcontents

\noindent
{\it Abstract.}
The object of this paper is the   theta divisor of the compactified Jacobian of a nodal curve.
We determine its irreducible components and
 give it a geometric interpretation. 
A characterization of  hyperelliptic  irreducible stable curves
  is appended as an application.

\


\

 
\section{Introduction}

Let $X$ be a connected, projective  curve of arithmetic genus $g$ and $\Pic^dX$ its degree-$d$ Picard variety,
parametrizing line bundles of degree $d$. If $X$ is smooth
$\Pic^dX$ is isomorphic to an abelian variety  and it is endowed with a principal polarization:
the theta divisor. If $d=g-1$ the theta divisor can be intrinsically defined as the locus
of $L\in \Pic^{g-1}X$ such that $h^0(X,L)\neq 0$.

If $X$ is singular,   
$\Pic^dX$  may fail to be projective, so one often needs to replace it with some projective
analogue, a so-called ``compactified jacobian", or ``compactified Picard variety". We shall
always assume that $X$ is reduced,  possibly reducible, and has at most nodes as singularities.

Although 
there exist several different constructions of compactified  jacobians in the literature, recent work of
V.Alexeev shows that in case $d=g-1$,
there  exists a ``canonical" one. More precisely, in \cite{alex}  the
 compactifications of T.Oda and C.S.Seshadri \cite{OS}, of C.Simpson \cite{simpson}, and
of \cite{caporaso}, are shown to be  isomorphic if $d=g-1$, to be endowed with an ample Cartier divisor,
the theta divisor $\Theta(X)$,
and 
to behave
consistently  with the degeneration theory of principally polarized abelian varieties. 

Some first results on the
theta divisor of the (non compactified)  generalized jacobian of any nodal curve were obtained  by A.Beauville;
see
\cite{beau}. Years later, A.Soucaris and E.Esteves independently constructed the theta divisor (as a Cartier, ample divisor)
on    the compactified jacobian of an irreducible curve; 
see \cite{soucaris} and  \cite{esttheta}.
The case of a  reducible, nodal curve was carried out in  
\cite{alex}.
As a result, today we know that, in degree $(g-1)$, 
the compactified  Picard variety
of any nodal curve  has a polarization,    the
theta divisor, 
such that the pair ({\it Compactified  Jacobian, Theta  Divisor} ) is a
semiabelic stable pair in the sense of \cite{alex1}. Furthermore,
 the above holds in the relative setting, i.e. for families of nodal curves.

These recent developements revive  interest in the theory of Brill-Noether varieties   for singular curves, 
of which 
 the theta divisor  is one of the principal
objects.

The purpose of this paper is to investigate  the geometry 
and the  modular meaning of $\Theta(X)$  more closely. Our  first result (Theorem~\ref{irr}) describes its irreducible
components, establishing that every irreducible component of the compactified jacobian
contains a unique irreducible component of the theta divisor, 
unless $X$ has some separating node (see
\ref{thetadef}); in particular, we characterize singular curves whose theta
divisor is irreducible (in \ref{irrth}). In more technical terms,
we prove  that for every fixed ``stable" multidegree
 (cf. Definition~\ref{baldef}) the theta divisor has a unique irreducible component.
This  result is sharp in the sense that irreducibility  fails for non stable multidegrees 
(see   Examples~\ref{ssred1}). The idea and the strategy of the proof are
described in \ref{idea}.

We prove the irreducibility Theorem~\ref{irr} using the Abel map,
namely, the rational map from $X^{g-1}$ to $\Pic^{g-1}X$, sending $(p_1,\ldots, p_{g-1})$ to 
$[\O_X(\sum p_i)]$.  
As a by-product,  the theta divisor is shown to be
the closure of the image of the Abel map, for every stable multidegree.
This fact, albeit   trivial for smooth curves, fails if the multidegree is not stable
(see Proposition~\ref{pureW} for a non semistable multidegree, 
and Example~\ref{ssred1} for a strictly semistable one).

In the second part of the paper we  
concentrate on the geometric interpretation of $\Theta(X)$ and precisely describe the objects it
parametrizes. 
In
 Theorem~\ref{Tstr} we   exhibit a  stratification   by means of the theta
divisors of the  partial normalizations of $X$.
We wish to  observe that very similar stratifications have been proved to exist for
several other compactified spaces, associated to singular curves (see \ref{Pstr}, or Theorem 7.9 in \cite{cner}, for
example). 
It is thus quite natural to ask whether all  compactified moduli spaces associated to a singular curve 
 admit   an analogous stratification, or whether some general rules governing such a phenomenon
exist. These questions are open at the moment.

Our stratification of $\Theta(X)$ yields a  description   in terms of effective  line bundles on the partial
normalizations of
$X$, or (which turns out to be the same) in terms of line bundles on semistable curves stably equivalent to $X$.

In the  final part, we apply our  techniques to
generalize to  singular curves the characterization of smooth hyperelliptic curves via
  the singular locus of their theta divisor; recall  that $\Theta(C)_{\text{sing}}=W^1_{g-1}(C)$
for every smooth curve $C$
of genus $g\geq 3$. Furthermore   $C$   is
hyperelliptic  if $\dim W^1_{g-1}(C)=g-3$, and non hyperelliptic if
$\dim W^1_{g-1}(C)=g-4$; we prove that the same  holds if $X$ is an
irreducible singular curve  (Theorem~\ref{W1thm}), but fails if $X$ is
reducible   (see \ref{nonhyp}). On the other hand the relation
between $\Theta(X)_{\text{sing}}$ and
$W^1_{g-1}(X)$  (and more  generally   $W^r_{g-1}(X)$)  i.e. a Riemann Singularity Theorem  for singular curves,  is not  known
 and it would be very interesting
to know it.

The paper consists of five sections. The first contains   preliminaries and basic definitions;
the second is mostly made  of  technical results. In the third section we prove the irreducibility theorem and
study the dimension of the image of the Abel map (Proposition~\ref{dimA}).
In the fourth section we describe the compactification of the theta divisor inside the compactified jacobian.
The fifth section  contains the application to singular hyperelliptic curves.

I wish to thank Juliana Coelho and the referee for several useful remarks. 
\subsection{Notation and Conventions}
\begin{nota}{}
\label{notX}
We work over an algebraically closed field $k$.
By the word ``curve" we  mean a reduced, projective curve over $k$.

Throughout the paper $X$ will be a  connected nodal curve of arithmetic genus 
$g$, having $\gamma$
irreducible components and $\delta$ nodes. We call
$\nu:Y\to X$  the normalization of $X$, so that $Y=\coprod_{i=1}^\gamma C_i$ with $C_i$  smooth of
genus
$g_i$, and
$X=\cup\overline{C_i}$ with $\overline{C_i}=\nu (C_i)$. Recall that
$
g=\sum_{i=1}^\g g_i +\dd -\g +1.
$
Observe that this formula holds regardless of $X$ being connected.

We denote by $\sing$ the set of nodes of $X$. For any set of  nodes of $X$,
 $S\subset \sing$,    set $\#S=\dS$ 
and $S=\{n_1,\ldots,n_{\dS}\}$.
The   normalization of $X$ at exactly the nodes in $S$ will  be denoted
$
\nS:\YS\la X
$
and $\gS$ will be the number of connected components of $\YS$; thus 
$
\YS=\coprod_1^{\gS}Y_i
$
 with $Y_i$ a connected curve of arithmetic genus $g_{Y_i}$.
We have
\begin{equation}
\label{genusS}
g=\sum_{i=1}^{\gS} g_{Y_i} +\dS -\gS +1 
\end{equation}
and, denoting $g_{\YS}=p_a(\YS)$
\begin{equation}
\label{genusY}
g_{\YS}=g-\dS =\sum_{i=1}^{\gS} g_{Y_i} -\gS+1.
\end{equation}
For every $j=1\ldots, \dS$  (or for every $n\in S$)  we set
\begin{equation}
\label{br}
\nS^{-1}(n_j)=\{q^j_1,q^j_2\}\  \  \  \  \  (\text{or }\  \nS^{-1}(n)=\{q_1,q_2\}).
\end{equation}
\end{nota}
\begin{nota}{}
\label{graph}
The dual graph of a nodal curve $Y$, denoted  $\Gamma_Y$, has    vertices  the irreducible
components of $Y$ and edges the nodes of $Y$. A node  
lying in a unique irreducible component $C_i$
is to a loop of $\Gamma_Y$ based at the vertex $C_i$; 
a node 
lying in $C_i\cap C_j$ is an edge joining the vertices $C_i$ and $C_j$.
\end{nota}
\begin{nota}{}
\label{notpic}
The degree-$d$ Picard variety $\Pic^dX$ has a decomposition into connected/ irreducible components:
$
\Pic^dX=\coprod_{\md\in \Z^\gamma:|\md|=d}\picX{\md}
$,
where $\picX{\md}$ is the variety of isomorphism classes of line bundles of multidegree $\md$.

Let $\nS:\YS \la X$ be as in \ref{notX}.
Consider the pull-back map
$$
\Pic X \stackrel{\nS^*}{\la } \Pic \YS\cong \prod_{i=1}^{\gS} \Pic Y_i \la 0.
$$
We shall usually identify $\Pic \YS\cong \prod \Pic Y_i$
without mentioning it. 

Let $M\in \Pic \YS$, then the fiber  over $M$ will be denoted
\begin{equation}
\label{FM}
\FM :=\{L\in \Pic X: \nS^*L=M\}\cong (k^*)^{\dS-\gS+1}.
\end{equation}
\end{nota}
\begin{nota}{}
\label{Lc}
We shall now describe the
isomorphism $\FM\cong (k^*)^{\dS-\gS+1}$ explicitely to fix some conventions.
Let us simplify the notation  by omitting the subscript $S$ (so, $\dd =\dS$, $Y=\YS$...).
Assume first that $Y$ is connected. 

Let $\mc = (c_1,\ldots, c_\dd)\in (k^*)^\dd$; $\mc$ determines a unique $L\in \Pic X$ such that
$\nu^*L=M$ as follows.
For every $j=1,\ldots ,\dd$ consider the two fibers of $M$ over $q^j_1$ and $q^j_2$ (recall that
$\nu(q^j_1)=\nu(q^j_2)=n_j)$, and fix an isomorphism between them. We define a line bundle $L=L^{(\mc)}$
on
$X$ which pulls back to $M$, by gluing
$M_{q^j_1}$ to $M_{q^j_2}$  via the isomorphism 
$$
M_{q^j_1}\stackrel{\cdot c_j}{\la} M_{q^j_2}
$$ 
given by multiplication by $c_j$.
Conversely, every $L\in \FM$ is of type $L^{(\mc)}$.

Now let $Y$ have $\gamma $ connected components; note that, since $X$ is connected,
 we always have $\gamma -1\leq
\dd$.
There exist some   subsets 
$T\subset S$  such that $\#T=\gamma -1$ and such that if we remove from $\Gamma _X$
every node that is not in $T$, the remaining graph 
is a connected tree
(a so-called {\it spanning tree} of $\Gamma _X$).

Let us fix one such $T$ and order the nodes in $S$ so that the last $\gamma -1$ are in $T$, 
i.e.
$
S=\{n_1,\ldots, n_{\dd}\}=\{n_1,\ldots, n_{\dd - \gamma +1}\}\cup T.
$
Now factor $\nu$ as follows
$$
\nu:Y\stackrel{\nu_{T}}{\la} Y'\stackrel{\nu '}{\la} X
$$
so that $\nu '$ is the partial normalization of $X$ at $S\smallsetminus T$ and $\nu _T$
the normalization at the nodes of $Y'$ preimages of the nodes in $T$.
For example, if $S=\sing$ (i.e. if $Y$ is smooth) then   $Y'$ is a curve of compact type.
The pull-back map $\nu_{T}^*$ induces an isomorphism
$
\Pic Y' \cong  \Pic Y,
$
i.e.  different gluing data determine isomorphic line bundles on $Y'$.

Now, to construct  the fiber of $\Pic X\to \Pic Y'$ over $M'$ we proceed as in the previous part.

Summarizing, for every $\mc \in (k^*)^{\dd  -\gamma +1}$ we associate a unique $L^{(\mc)}\in \Pic Y$;
since the gluing data over the nodes in $T$ is irrelevant, we shall fix $c_j=1$ 
if $j\geq \dd -\gamma$ and use that as gluing constant over $T$.

Finally
observe that a section $s\in H^0(Y,M)$ descends to a section
$\overline{s}\in H^0(X,L^{(\mc)})$ if and only if for every $j=1,\ldots ,\dd$ we have
\begin{equation}
\label{glue}
s(q^j_2)=c_js(q^j_1).
\end{equation}
\end{nota}

\subsection{Brill-Noether varieties and Abel maps}
\begin{nota}{}
\label
We recall some basic facts about Brill-Noether varieties for smooth curves, following
the notation of
\cite{ACGH} to which we refer for  details.

Let $C$ be a smooth connected curve of genus $g\geq 0$, and let $d$ and $r$ be nonnegative integers.
The set $
\W _d^r(C):=\{ L\in \Pic^dC: h^0(C,L)\geq r+1\}
$
has an algebraic structure  and  is called a Brill-Noether variety. 
It is closely related to  the Abel map in degree $d$ of $C$, that is the map
\begin{equation}
\label{absm}
\begin{array}{lccr}
\alpha^d_C: & C^d&\la &\Pic^d(C)\\
&(p_1,\ldots,p_d) &\mapsto &\O_C(\sum_1^dp_i).
\end{array}
\end{equation}
Then $\im \alpha^d_C\subseteq \W _d^0(C)$ for all $d\geq 0$
(see \ref{r(d)} for when equality occurs).
Note that $\W _d^r(C)$ may fail to be irreducible, so when talking about its
dimension we will mean the maximum dimension of its components.
The following is well known(\cite{ACGH} Lemma 3.3 Ch.IV) 
\begin{fact} 
\label{rho}
If $r\geq d-g$,   every irreducible component of $\W _d^r(C)$ has dimension at least equal to
$\BN:=g-(r+1)(r-d+g)
$.
If $r\leq d-g$ then 
%(by Riemann-Roch) 
$\W _d^r(C)=\W _d^{d-g}(C)$.
\end{fact}

There is also a simple upper bound 
\begin{equation}
\label{dimW}
\dim \W _d^r(C)\leq\min \{d-r,g\}.
\end{equation}
Indeed, if $d-r\leq g$, it suffices to look at the
Abel map of degree $d$ to obtain that $\dim \W _d^r(C)\leq d-r$ 
 (cf. \cite{ACGH} Prop. 3.4 Ch.IV).
If $d-r\geq g$ then, by Riemann-Roch,   $\dim \W _d^r(C)=g$ 

\begin{remark}
\label{r(d)}
Denote by $r(d) $ the dimension of a general (non empty)
complete linear system of degree $d$. i.e.
if $d\leq g$ set $r(d)=0$, if $d\geq g$ set  $r(d)=d-g$.
Note that $\W_d^{r(d)}(C) = \im\alpha^d_C$.
Now,  
$
\min \{d-r(d),g\}=\min\{d,g\}
$
and 
\begin{displaymath}
\dim \W _d^r(C)\left\{ \begin{array}{ll}
=\min\{d,g\} &\text{ if } r\leq r(d)\\
< \min\{d,g\} &\text{ if } r> r(d).\\
\end{array}\right.
\end{displaymath}
To see that, assume first that $r\leq r(d)$, then $\W _d^r(C)=\W _d^{r(d)}(C)$
 by Riemann-Roch, so we may assume that $r=r(d)$.
Now computing gives
$
\rho(d,g,r(d))=\min\{d,g\}
$,
so by fact~\ref{rho} and (\ref{dimW}) we get $\dim \W _d^r(C)=\min\{d,g\}$.
The  case $r> r(d)$  follows from (\ref{dimW}) and   the fact that
$\min\{d-r,g\}<\min\{ d-r(d),g\}$.
\end{remark}
\end{nota}
\begin{nota}{}
\label{notW}
For a nodal curve $X$ of genus $g$ having $\gamma$ irreducible components, for any 
$\md \in \Z^{\gamma}$ and $r\geq 0$, we set
$
\Wmdr=\{L\in \picX{\md}: h^0(X,L)\geq r+1\}
$ 
and for any $d\in \Z$,
$
W_d^r(X):=\coprod_{|\md|=d}\Wmdr
$.
In case $r=0$  the superscript $r=0$ is usually omitted.
In particular 
$$
\W_{g-1}(X):=\{L\in \picX{g-1}: h^0(X,L)\geq 1\}=\coprod_{|\md|=g-1}\Wmd.
$$
With the notation of \ref{notpic}, if $\nS:\YS\to X$ is a partial normalization and $M\in \Pic\YS$, the fiber of $W_d^r(X)$ over $M$  will be denoted  (recall (\ref{FM}))
\begin{equation}
\label{WMr}
\WMr:=\{L\in \FM: h^0(X,L)\geq r+1\}
\end{equation}
and $
\WM:=\{L\in \FM: h^0(X,L)\geq 1\}.
$
\begin{remark}
\label{sconn} 
The above definitions make sense  also for non connected curves.
Consider a disconnected  curve, $Y=\coprod _{i=1}^\gamma C_i$ where $C_i$ is smooth and connected
(or more generally $C_i$  irreducible) of genus $g_i$.
For any $\md \in \Z^\gamma$,
the variety $\W _{\md}(Y)$ is easily described in terms of the $C_i$:

\begin{displaymath}
\W _{\md}(Y) = \left\{ \begin{array}{ll}
\prod_{i=1}^\gamma \Pic^{d_i}C_i &\text{ if } \exists i:\  d_i\geq g_i\\
\bigcup_{j=1}^\gamma  \Bigr(\W _{d_j}(C_j) \times \prod_{i\neq j, i=1,\ldots, \gamma}
\Pic^{d_i}C_i\Bigl) &\text{ if }  \forall i:\  d_i\leq g_i-1.\\
\end{array}\right.
\end{displaymath}
\end{remark}

We shall need the following very simple
\begin{lemma}
\label{nbp} Let $S\subset \sing $,
$\nS:\YS \to X$  the normalization of $X$ at $S$ and $p\in X\smallsetminus S$.
Let $M\in \Pic \YS$ and assume that $M$ has no base point in $\nS^{-1}(S\cup p)$.
Then there exists $L\in \WM$ such that $L$ has  no base point in $p$.
In particular, if $M$ has no base point over $S$ then $\WM$ is nonempty.
\end{lemma}
\begin{proof}
To say that $M$ has no base point in $\nS^{-1}(S\cup p)$ is to say that there exists $s\in H^0(\YS,M)$ such that
$s(q)\neq 0$ for every $q\in \nS^{-1}(S\cup p)$.
We can use $s$ to construct a line bundle $L\in \WM$ by identifying the two fibers over pairs of corresponding
branches. More precisely, with the notation of \ref{Lc} (\ref{glue})
for every $n_j\in S$ call $q_1^j,q_2^j$ the branches over $n_j$. Then set
$c_j:=s(q^j_2)/s(q^j_1)$ and define $L=L^{(\mc)}$.
It is clear that $s$ descends to a nonzero section $\ov{s}$ of $L$ and that $\ov{s}(p)\neq 0$.
\end{proof}
\end{nota}

\begin{nota}{\it Abel maps.}
\label{abel}
We now introduce the Abel maps of a singular curve.
Recall (see \ref{notX}) that $X=\ov{C}_1\cup\ldots \cup \ov{C}_{\gamma}$ denotes the decomposition of $X$ into 
irreducible components. For every $\md =(d_1,\ldots, d_{\gamma})$ such that $d_i \geq 0$ we set
$\Xmd=\ov{C}_1^{d_1}\times\ldots \times \ov{C}_{\gamma}^{d_{\gamma}}$.
Now denote $\dfib=X\smallsetminus \sing$
the smooth locus of $X$.
The normalization map 
$
Y=\cup C_i\stackrel{\nu}{\la} X=\cup \ov{C}_i
$ 
induces an isomorphism of $\dfib$ with $Y\smallsetminus \nu^{-1}(\sing)$. We shall identify 
$\dfib = Y\smallsetminus \nu^{-1}(\sing)$ and  denote
$
\dot{C_i}: = \ov{C}_i\cap \dfib .
$
Finally, set
$$
\dXmd:=\dot{C_1}^{d_1}\times\ldots \times \dot{C_{\gamma}}^{d_{\gamma}}\subset \Xmd
$$
so that $\dXmd$ is a smooth irreducible variety of dimension $|\md|$, open and dense in $\Xmd$.
Set $d=|d|$, then we have a regular map

\begin{equation}
%\label{}
\begin{array}{lccr}
\amd: & \dXmd&\la &\picX{\md} \\
&(p_1,\ldots,p_d) &\mapsto &\O_X(\sum_1^d p_i)
\end{array}
\end{equation}
which we call the {\it Abel map of multidegree $\md$}. 
 We denote
$$
\Amd:=\overline{\amd(\dXmd)} \subset \picX{\md}.
$$

\begin{lemma}
\label{Ad} Let $X$ be a (connected, nodal) curve of genus $g\geq 0$.
For every $d\geq 1$ and every multidegree $\md$ on $X$ such that
$\md \geq 0$ and $|\md |=d$ we have
\begin{enumerate}[(i)]
\item
\label{Ad1}
$\Amd$ is  irreducible  and  $\dim \Amd\leq \min\{d,g\}$;
\item
\label{Ad2}
$\Amd\subset \Wmd$.
\end{enumerate}
\end{lemma}
\begin{proof} Obvious.
%(\ref{Ad1}) is obvious, as $\Xmd$ is irreducible of dimension $d$ and $\dim \picX{\md}=g$.
%
%For any $(p_1,\ldots,p_d)\in \dXmd$ we  have  $h^0(X, \sum p_i)\geq 1$, hence
%(\ref{Ad2}) follows  by uppersemicontinuity of $h^0$.
\end{proof}

We shall see  that  strict inequality in (\ref{Ad1}) does occur (cf. Proposition~\ref{dimA}).

\subsection{Stability and semistability}
As we said in the introduction, there exist various  modular descriptions for a compactified Picard variety,
and they are   equivalent if $d=g-1$. 
We shall give the  complete description later, in \ref{descgen}. For now it is enough to recall that,
for every nodal curve $X$,
 the compactified Picard variety in degree $g-1$,  $\PXgb$, is   a union of 
 (finitely many) irreducible  $g$-dimensional components 
each of which contains as an   open subset a copy of the generalized Jacobian of $X$.
 To study the irreducible components of the theta divisor of $\PXgb$ 
 there is no need to consider its boundary points.
 This explains why we chose to postpone the complete description of $\PXgb$; see \ref{descgen}.

 So,  now  only the open smooth
locus of $\PXgb$
will be described, using line bundles of  ``stable" multidegree
 on the normalization of $X$ at its separating nodes. 

There exist two different,
equivalent definitions
of semistability and stability (\ref{baldef} and \ref {beau} below); the simultaneous use of the two
is a good tool to overcome technical  difficulties of  combinatorial type.
\end{nota}
\begin{nota}{\it Stability: Definition 1.}
\label{baldef}
Let $Y$ be a  nodal curve of arithmetic genus $g$ having $\gamma$ irreducible components. Let
$\md\in
\Z^{\gamma}$ be such that
$
|\md|=g-1.
$ 
\begin{enumerate}[(a)]
\item
\label{bal}
We call $\md$   {\it semistable}
%
%\footnote{{\it Semistable} is the same as {\it balanced} in \cite{cner}.}
%
if for every subcurve (equivalently, every connected subcurve)
$Z\subset Y$ of arithmetic genus $g_Z$
we have
\begin{equation}
\label{BIm}
d_Z\geq g_Z -1
\end{equation}
where $d_Z:=|\md_Z|$.
The set of   semistable multidegrees  on   $Y$ is denoted
$$
\BYs :=\{\md\in \Z^{\gamma}:\   |\md|=g-1,\  \md\   \text{ is semistable}\}.
$$
\item
\label{stbal}
Assume $Y$ connected.
If   $Y$ is irreducible, or  if strict inequality holds in (\ref{BIm})
for every (connected) subcurve
$Z\subsetneq Y$, then  $\md$ is called  {\it stable}.
%
%\footnote{{\it Stable} here is not the same as {\it stably balanced} in \cite{cner}}
%
If $Y$ is not connected, we say that $\md$ is stable if its restriction to every connected
component of $Y$ is stable.
We denote
$$
\BY :=\{\md\in \Z^{\gamma}:\   |\md|=g-1,\  \md\   \text{ is stable}\}\subset \BYs .
$$
\end{enumerate}

\end{nota}

\noindent
We shall also use the following equivalent definition,
originating from \cite{beau}.
\begin{nota}{\it Stability: Definition 2.}
\label{beau} 
Fix $Y$ and $\md$ as in \ref{baldef}.
\begin{enumerate}[{\it (A)}]
\item
\label{beaubal}
 $\md$   is {\it semistable}
if the dual graph $\Gamma_Y$ of $Y$ (cf. \ref{graph}) can be oriented in such a way that, denoting by $b_i$
the number of edges pointing at the vertex corresponding to the irreducible component $C_i$
of $Y$, then
$$
d_i=g_i-1+b_i
$$
where $g_i$ is the geometric genus of $C_i$ (so that $g_i=p_a(C_i)-\#(C_i)_{\text{sing}}$).
\item
\label{beaust}
Assume $Y$ connected. Then  $\md$ is stable if $\Gamma_Y$ admits an orientation satisfying
(\ref{beaubal}) and such that there exists no proper
subcurve
$Z\subsetneq Y$ such that  the edges   between $\Gamma_Z$ and $\Gamma _{Z^C}$ go all in
the same  direction, where  $Z^C:=\overline{Y\smallsetminus Z}$.
\end{enumerate}

The equivalence of definitions \ref{beau} and \ref{baldef} is Prop.3.6 in \cite{alex}.
The version  given in \ref{beau} (\ref{beaubal}) is due to A.Beauville, who used it in
\cite{beau} to
define and study the theta divisor of a  generalized jacobian
(In \cite{beau} Lemma (2.1) the dual graph  is  without loops by definition,
whereas we need to include loops. This
explains the  difference between our definition and that of \cite{beau}).

Version~\ref{baldef} actually extends to all degrees (other than degree $g-1$);
 it originates from D.Gieseker's construction of $\mgbar$ and
is crucial in  \cite{caporaso}
(where (\ref{BIm}) is generalized by the so-called ``Basic Inequality").
V.Alexeev proved   that the Basic  Inequality 
 yields
the modular description of the compactified jacobians constructed 
by Oda-Seshadri
and by C.Simpson using  different approaches (see \cite{alex} 1.7 (5)). 
More details about this definition and its connection with Geometric Invariant Theory
will be given in Section~\ref{comp}.

\begin{remark}
\label{balrk} 
\begin{enumerate}[(i)]
\item
\label{}
Applying inequality (\ref{BIm}) to all subcurves, we get that $\md$ is semistable if and only if for
every connected $Z\subset Y$
\begin{equation}
\label{BI}
p_a(Z) -1 \leq d_Z \leq p_a(Z) -1 +\#Z\cap Z^C.
\end{equation}
If $X$ is connected, $\md$ is stable if and only if strict inequalities hold in (\ref{BI}) for all $Z$.
\item
\label{}
If $\md\in \BXs$ 
and $V\subset X$ is a subcurve such that $d_V=g_V-1$,
then $\md_V$ is semistable on $V$.
\item
\label{}
If $\md$ is stable, then $\md \geq 0$.
\end{enumerate}\end{remark}

\begin{remark}
\label{he} The following convention turns out to be useful. Given a graph $\Gamma $ (e.g. $\Gamma
=\Gamma _Y$), every edge $n$ determines two half-edges, 
denoted $q_1^n$ and $ q_2^n$  (corresponding to the two branches of the node $n$ of $Y$). If $\Gamma $
is oriented we call $q_1^n$ the starting half-edge  of $n$ and
$q_2^n$ the ending one.
 \end{remark}

$\BXs$ is never empty (by \cite{cner} Prop. 4.12). On the other hand we have

\begin{lemma}
\label{empty} 
%%Let $X$ be a connected  curve. 
$\BX=\emptyset$ if and only if $X$ has a separating
node.
\end{lemma}
\begin{proof}
If $X$ has a separating node,   $n$, then $X=X_1\cup X_2$ with $X_1\cap X_2=\{n\}$.
Let $\md \in \BXs$, using (\ref{BI}) we have
$ 
p_a(X_i)-1\leq d_{X_i} \leq p_a(X_i),
$ 
so that  strict inequalities cannot simultaneously occur. Hence $\md$ is not stable.

Conversely, assume that $X$ has no separating node.
We shall  use  Definition~\ref{beau}, and prove that 
the dual graph of $X$,  $\Gamma=\Gamma _X$, admits a ``stable orientation"
(i.e. an orientation satisfying  (\ref{beaust})).
We use induction on the number $\dd$ of nodes that lie in two different irreducible components
(the only nodes that matter), i.e. induction on the number of edges that are not loops.
If $\dd=1$ there is nothing to prove (the edge is necessarily separating),
if $\dd=2$ then $\Gamma$ has two vertices so the statement is clear.

Let $\dd \geq 2$, pick an edge $n$ and let $\Gamma '=\Gamma -n$; thus $\Gamma '$ is
connected. If $\Gamma '$ has no separating edge, by induction $\Gamma '$ admits a stable
orientation, hence so does $\Gamma$, of course.
Denote $n_1,\ldots, n_t$ the separating edges of
$\Gamma '$. The graph
$$
\Gamma '-\{n_1,\ldots, n_t\}= \Gamma  -\{n_0,n_1,\ldots, n_t\}
$$
where $n=n_0$, has $t+1$ connected components, $\ov{\Gamma _0}$,\ldots, $\ov{\Gamma _t}$,
each of which is free from separating edges.

We claim that the image 
${\Gamma _i}\subset \Gamma$ of each
$\ov{\Gamma _i}$  contains exactly two of the edges $n_0,n_1,\ldots, n_t$.

Indeed, if (say) $\Gamma_1$ contains only one $n_i$ with $i\geq 1$, call it $n_1$
and call $\Gamma_2$ the other $\Gamma_i$ containing $n_1$. Then $n_0$ connects $\Gamma_1$ with
$\Gamma_2$
(for otherwise $n_1$ would be a separating node of $\Gamma$ which is not possible). Hence $\Gamma_1$
contains $n_0$ and $n_1$.

If   $\Gamma_1$ contains two $n_i$ with $i\geq 1$, call them $n_1$ and $n_2$, let $\Gamma_2$ and $\Gamma
_3$ be such that
$n_i\in \Gamma_1\cap\Gamma_{i+1}$, $i=1,2$. Then $n_0$ connects $\Gamma_2$ and $\Gamma _3$,
thus $n_0\not\in\Gamma_1$. 
Therefore $\Gamma_1$ contains only $n_1$ and $n_2$.

If   $\Gamma_1$ contains three $n_i$, $i\geq 1$, call them $n_1$, $n_2$ and $n_3$, let $\Gamma_2$,
$\Gamma _3$ and
$\Gamma _4$ be such that
$n_i\in \Gamma_1\cap\Gamma_{i+1}$. Now $n_0$ is contained in at most two $\Gamma_i$,
so 
say $n_0\notin\Gamma_4$ (say), but then $n_3$ is a separating node of $\Gamma$, which is a
contradiction.
Therefore, up to reordering the $\Gamma _i$, we can assume that
$$
n_i\in  {\Gamma _i}\cap\Gamma_{i-1},\  \  i=1,\ldots, t,t+1=0.
$$
We now define an orientation on $\Gamma$ by combining the stable orientation on each
$\Gamma_i$ with  each edge $n_i$ oriented from $\Gamma_{i-1}$ to $\Gamma _i$.
It suffices  to prove that this is a stable orientation on $\Gamma$.

Indeed: let $Z\subset X$ and $\Gamma _Z\subset \Gamma$ the corresponding graph.
If for some $i$ we have $\emptyset\neq \Gamma_Z\cap \Gamma _i\subsetneq\Gamma_i$,
then inside $\Gamma _i$ there are edges  both starting from and ending in $\Gamma _Z$.
So the same holds in $\Gamma$ and we are done.
Hence we can assume that for every $i$ either $\Gamma_i\subset \Gamma _Z$ or
$\Gamma _Z\cap \Gamma _i=\emptyset$. Therefore
$$
\Gamma_Z\cap \Gamma _{Z^C}\subset \{n_0,n_1,\ldots, n_t\}.
$$
We can thus reduce ourselves to consider the 
 graph obtained by contracting every $\Gamma _i$ to a point. This is of course a cyclic graph
with $t+1$ vertices and $t+1$ edges $\{n_0,n_1,\ldots, n_t\}$, oriented cyclically.
This is a stable orientation, so we are done.
\end{proof}

\begin{example}
\label{Ssplit}
Let $X$ be a nodal connected curve of genus $g$, $\sep\subset\sing$ the set of its separating nodes
and $\widetilde{X}\to X$ the normalization of $X$ at $\sep$.
Assume $\#\sep=c-1\geq 1$ so that $\widetilde{X}$ has $c$ connected components $X_1,\ldots, X_c$
and $X_i$ is free from separating nodes for every $i=1,\ldots,c$. Thus $\Sigma(X_i)\neq \emptyset$ and
$$
\Sigma(\widetilde{X})=\Sigma(X_1)\times\cdots \times \Sigma(X_c).
$$
Indeed, set $g_i:=p_a(X_i)$, then
$p_a(\widetilde{X})-1=(g-c+1)-1=\sum_{i=1}^c(g_i-1)$,  and $\md\in \Sigma(\widetilde{X})$
if and only the restriction  of $\md$ to $X_i$ is stable on $X_i$.
\end{example}

\begin{prop}[Beauville]
\label{pureW}  
Let $X$ be a (connected, nodal) curve of genus $g\geq
1$, and let $\md\in \Z^{\gamma}$ be such that $|\md|=g-1$.
\begin{enumerate}[(i)]
\item
\label{}
 $\md$ is semistable iff there exists $L\in \picX{\md}$ such that
$h^0(X,L)=0$.
\item
\label{}
If $\md$ is semistable then every irreducible component of $\Wmd$ has dimension $g-1$.  
\item
If $\md$ is not semistable then  $\Wmd=\picX{\md}$. 
\end{enumerate}
\end{prop}
See Lemma (2.1) and Proposition 2.2 in \cite{beau}.
\end{nota}

\begin{nota}{}
\label{idea} Our first theorem (Theorem~\ref{irr}) states that, if $\md$ is stable, then $\Wmd$ is irreducible
and equal to $\Amd$. The proof's 
 strategy  is the following. We know, by the above Proposition~\ref{pureW},
that every irreducible component of $\Wmd$ has dimension $g-1$; we also know that $\Amd$ is irreducible.
We shall prove that if $W$ is an irreducible component of $\Wmd$, not dominated by the image of the Abel map,
then $\dim W\leq g-2$, and hence $W$ must be empty.

To do that we consider the normalization $\nu:Y\to X$ and the pull back map:
$\nu^*:\Pic X\to \Pic Y$. The dimension of $W$ is then
studied by fibering $W$ using  $\nu^*$, and bounding the dimensions of the image and the fibers.

An important    point is to show that, on the one hand,
the divisors on $Y$ supported  over the nodes of $X$
 impose independent conditions on the general line bundle $M\in \Pic^{\md}Y$; see  Lemma~\ref{icst}.
On the other hand, if $M\in \Pic Y$ 
has this property (i.e
 divisors   supported  in $\nu^{-1}(\sing)$
 impose independent conditions on it), then  the dimension of the locus of 
$L\in W_M(X)$ which do not lie in the image of the Abel map 
is small, hence the dimension of the fiber of $W$ over $M$
is small;
see Proposition~\ref{nvr} and Corollary~\ref{WSgen}.
\end{nota}

\section{Technical groundwork}

\subsection{Basic estimates}

Recall the set-up of \ref{notX}.

\begin{prop}
\label{GS} Fix $\nS:\YS\la X$  and let
 $M\in \Pic\YS$.
\begin{enumerate}[(i)]
\item
\label{}
For every $L\in \Pic X$ such that $\nS^*L=M$ we have
\begin{equation}
\label{eff}
h^0(\YS,M)-\dS\leq h^0(X,L)\leq h^0(\YS,M).
\end{equation}

 \item
\label{GSW}
Let  $h^0(\YS,M)\geq \dS$. Assume 
that  for some  $h:\{1,\ldots,\dS\}\to \{1,2\}$,
\begin{equation}
\label{ind}
h^0(\YS, M(-\sum_{j=1}^{\dS} q^j_{h(j)}))=h^0(M)-\dS.
\end{equation}
Then $\WM$ is of pure dimension:
\begin{displaymath}
\dim \WM=\left\{ \begin{array}{ll}
\dS - \gS &\text{ if } h^0(M)=\dS\\
\dS - \gS+1 &\text{ if } h^0(M)\geq \dS+1.\\
\end{array}\right.
\end{displaymath}
Moreover, the general element $L\in \WM$ satisfies
\begin{equation}
\label{GSh}
h^0(X,L) = \max \{h^0(\YS,M)-\dS, 1 \}.
\end{equation}
\end{enumerate}
\end{prop}

\begin{proof}
Throughout the proof we shall simplify the notation by omitting the index $S$, i.e. set $Y=\YS$,
$\dd =\dS$, $\nu =\nS$ and $\g=\gS$.

Let $L\in \FM$, then we have the exact sequence 
\begin{equation}
\label{S1}
0\la L \la \nu_*M\la \sum_{n\in S}
k_n\la 0
\end{equation}
and the associated long cohomology sequence
\begin{equation}
\label{C1}
0\to H^0(X,L) \stackrel{\alpha}{\la} H^0(Y,M)\stackrel{\beta}{\la} k^{\dd}\to H^1(X,L) \to H^1(Y,M)\to 0
\end{equation}
from which we immediately get the upper bound on $h^0(X,L)$ stated in (\ref{eff}).

Fix $M\in \Pic Y$, recall the description
of the fiber of $\nu^*$ over $M$ 
given in \ref{Lc}. 
Thus every $L\in \FM$ is of the form $L=L^{(\mc)}$ for some $\mc \in (k^*)^{\dd-\gamma +1}$.
For convenience, we  use the same  set-up of \ref{Lc}, in particular
we set $c_j=1$ for $\dd -\gamma +2\leq j \leq
\dd$.

To compute $H^0(X,L)$, set $
l=h^0(Y,M)
$ and  pick a basis $s_1,\ldots, s_l$ for $H^0(Y,M)$.
 Let $s\in H^0(Y,M)$, so $s=\sum_1^lx_is_i$
where $x_i\in k$. Now $s$ descends to a section of $L$ 
(i.e. $s$ lies in the image of $\alpha$ in (\ref{C1})) if and only if
\begin{equation}
\label{des}
\sum_{i=1}^lx_i\bigr( s_i(q^j_2)-c_js_i(q^j_1)\bigr)=0\hskip1.in \forall j=1,\ldots,\dd .
\end{equation}
The above is a linear system of $\dd$ homogeneous  equations in the $l$ unknowns $x_1,\ldots, x_l$.
The space of its solutions, $\Lambda(\mc)$, is  identified with $H^0(X,L^{(\mc)})$.
Now,
$\Lambda(\mc)$ is a linear subspace of $H^0(Y,M)$ of dimension at least $l-\delta$.
Hence $h^0(X,L)=\dim \Lambda(\mc)\geq l-\delta$, proving   (\ref{eff}).

Part (\ref{GSW}). Assume  $l=h^0(Y,M)\geq \dd$;
denote by
$A(\mc)$ the $\dd \times l$ matrix of the system (\ref{des}).
By what we said
\begin{equation}
\label{hdet}
h^0(X,L^{(\mc)})=\dim \Lambda(\mc)=l-rankA(\mc)
\end{equation}
and
\begin{equation}
\label{Wdet}
\WM\cong \{\mc : \Lambda(\mc )\neq 0\}=\{\mc: rank A(\mc) \lneq l\}.
\end{equation}

We shall  prove that $A(\mc)$ has rank $\dd$ unless $\mc$ lies in a proper closed subset
of
$(k^*)^\dd$. For that, we apply the assumption
(\ref{ind}) to
choose the basis for $H^0(Y,M)$ as follows. First, up to renaming each pair of branches we can assume that
$h(j)=1$ for every $j$. By (\ref{ind}) we can
pick  $\dd$ linearly independent $s_1,\ldots s_\dd\in H^0(M)$ such  that 
\begin{displaymath}
s_i(q^j_1)=\left\{ \begin{array}{ll}
1 &\text{ if } \  i=j\\
0 &\text{ if } j\neq i,\  (j=1,\ldots, \dd).\\
\end{array}\right.
\end{displaymath}

If $l>\dd$ we choose the remaining  basis elements however we like.
Set
$b^j_i:=s_i(q^j_2)\in k .$
Then the matrix $A(\mc)$ contains a $\dd\times \dd$ minor, $B(\mc)$,
(the minor given by the first $\dd$ columns)
whose diagonal is
$ 
(c_1-b^1_1,\ldots, c_\dd-b^\dd_\dd)
$,
and such that the $c_j$ do not appear anywhere else in $B(\mc)$.
Therefore the determinant of $B(\mc)$ is a nonzero polynomial in the $c_j$. This proves that the locus 
where the matrix has maximal rank (equal to $\dd$) is open, non empty.

Suppose $\dd =l$, then $B(\mc)=A(\mc)$. By (\ref{Wdet})
$\WM$ is naturally identified with the locus of points of $\FM$
where $\det A(\mc)$ vanishes.
We conclude that $\WM$ has pure dimension $\dim \WM=\dd-\gamma$ proving (\ref{GSW}).

Moreover, 
for a general $L^{(\mc)}\in \WM$, the rank of $A(\mc)$ is equal to $\dd -1$.
Indeed, by 
(\ref{Wdet}) $\WM$  is identified to the hypersurface,  $W$, of $k^\dd$
 where $\det A(\mc)$ vanishes. Call $A^i_j(\mc)$ the minor of $A(\mc)$
obtained by removing the $i$-th row and the $j$-th column, and set
$ 
U^i_j=\{\mc\in k^\dd: \det A^i_j(\mc)\neq 0\}$. We must prove that $W\cap U^i_j\neq \emptyset$
for some $1\leq i,j\leq \dd$. Suppose $c_1$ appears in $\det A(\mc)$. On the other hand
$c_1$ does not appear in $\det A^1_1(\mc)$, as 
$A^1_1(\mc)$ does not contain  $c_1$.
Hence $W\cap U^1_1\neq \emptyset$.

Therefore
 by (\ref{hdet})  we get
$h^0(X,L)=1$ proving (\ref{GSh}).

If $l>\dd$, then $\WM=\FM$ by part (\ref{eff}). Furthermore, by (\ref{hdet}) 
$$
h^0(X,L^{(\mc)} )=l-rank A(\mc)\geq l-\dd.
$$ 
By looking at the matrix $A(\mc)$, we see that  $h^0(X,L^{(\mc)})= l-\dd$ holds on the non empty open  subset
where
$\det B(\mc)$ does not vanish; this proves (\ref{GSh}).
\end{proof}

\begin{lemma}
\label{GSb} Let $\nu :Y\to X$ be the normalization of $X$  and let $\md\in \BXs$. For a general $M\in
\Pic^{\md}Y$ we have

\begin{enumerate}[(i)]
\item
\label{GSbM}
$h^0(Y,M)=\dd$;
\item
\label{indb}
$M$ satisfies condition (\ref{ind})  w.r.t. a suitable  $h:\{1,\ldots,\dd\}\to \{1,2\}$;
\item
\label{GSbW}
$\dim \WM = \delta - \gamma$;
\item
\label{GSbh} The general $L$ in $\WM$ satisfies $h^0(X,L)=1$.
\end{enumerate}
%(where $\delta=\#\sing$ and $\gamma$ is the number of irreducible components of $X$).
\end{lemma}
\begin{proof}
Using the notation of \ref{notX}, $Y=\coprod C_i$ with $C_i$   smooth of genus $g_i$, and
$X=\cup\overline{C_i}$. The fact that $\md$ is  semistable implies that $d_i\geq p_a(\overline{C_i})-1\geq g_i-1$
for every
$i=1,\ldots,\gamma$. Therefore for $M$ general in $\Pic^{\md}Y$
$$
h^0(Y,M)=\sum_i(d_i-g_i+1)=g-1-\sum_ig_i+\gamma=\delta.
$$
Let us prove (\ref{indb}).
We use definition \ref{beau} (\ref{beaubal}) of  a semistable multidegree;
 $\Gamma_X$ of $X$ can be oriented   so that, if $b_i$ denotes the number
of edges pointing at $C_i$, then for all $i=1,\ldots ,\gamma$
\begin{equation}
\label{BBI}
d_i=g_i-1+b_i .
\end{equation}
Any such orientation gives us a choice of branches over each node. Namely,
for every $n_j\in \sing$ we denote $q_2^j$ the branch corresponding to the ending half-edge of  $n_j$.
We claim that (\ref{ind}) holds with respect to the map   $h(j)=2$ for every $j$.
Indeed
$$
h^0(Y,M(-\sum _{j=1}^{\dd} q^j_2))=\sum_{i=1}^{\gamma} h^0(C_i,M(-\sum _{j=1}^{\dd} q^j_2)_{|C_i}).
$$
Now by (\ref{BBI})
\begin{equation}
\label{degM-E}
\deg_{C_i} M(-\sum _{j=1}^{\dd} q^j_2)=d_i-b_i=g_i-1
\end{equation}
hence ($M$ being general) 
$h^0(C_i,M(-\sum _{j=1}^{\dd} q^j_2)_{|C_i})=0$ for every $i$. We conclude that, by part (\ref{GSbM}),
$$
h^0(Y,M(-\sum _{j=1}^{\dd} q^j_2))=0=h^0(Y,M) -\dd
$$
so that (\ref{ind}) is satisfied. 
Now, applying \ref{GS}(\ref{GSW}),
we get $\dim\WM=\delta - \gamma$
and    $h^0(X,L)=1$ for a general $L\in \WM$.
So (\ref{GSbW}) and (\ref{GSbh}) are proved.
\end{proof}

\begin{cor}
\label{GSc} 
Let $\md\in \BXs$ and let 
 $L$ be a general line bundle in $\picX{\md}$. For every subcurve $Z\subseteq X$ we have
$
h^0(Z,L_Z)=d_Z-g_Z+1.
$
\end{cor}
\begin{proof}
It  suffices to assume $Z$ connected (by (\ref{genusY})).
Consider the normalization  $\nu :Y=\cup C_i\to X$ of $X$  and 
$\nu^*L=M=\Lv$ with $L_i\in \Pic^{d_i}C_i$.  Then
$L_i$ is general in $\Pic^{d_i}C_i$ (as $L$ is general in $\picX{\md}$); since $d_i\geq g_i -1$ 
(as $\md$ is semistable) we get that
$h^0(C_i, L_i)=d_i-g_i+1$.
Now, denote $Z^{\nu}\to Z$ the normalization of $Z$, 
order the irreducible components of $X$ so that the first $\gamma _Z$ are the irreducible components of $Z$,
set 
$S=Z_{\text{sing}}$, so that 
$g_Z=\sum_{i=1}^{\gamma _Z} g_i+\dS - \gamma _Z+1.$
Let
$M_{Z^{\nu}}$ be the restriction of $M$ to  $Z^{\nu}$, then
$$
h^0(Z^{\nu},M_{Z^{\nu}})=\sum_{i=1}^{\gamma_Z}h^0(C_i,L_i)=\sum_{i=1}^{\gamma_Z}(d_i-g_i+1)=d_Z-g_Z+\dS +1.
$$ 
Now, since $\md$ is semistable, $d_Z\geq g_Z-1$ hence $h^0(Z^{\nu},M_{Z^{\nu}})\geq \dS$.
Moreover, recall that by \ref{GSb}  (\ref{indb}) $M$ satisfies condition (\ref{ind}); 
it is straightforward to check
that the analogue holds for $M_{Z^{\nu}}$, i.e. for a suitable choice of branches,
$$
h^0(Z^{\nu}, M_{Z^{\nu}}(-\sum_{j=1}^{\dS} q^j_{h(j)}))=h^0(M_{Z^{\nu}})-\dS=0.
$$
This enables us to apply \ref{GS}(\ref{GSh})
 to $Z^{\nu}\to Z$,
 thus getting 
$$
h^0(Z,L_Z)=h^0(Z^{\nu},M_{Z^{\nu}})-\dS=d_Z-g_Z+\dS +1-\dS=d_Z-g_Z +1.
$$
\end{proof}

\subsection{Basic cases}
Recall the notation of \ref{notX}, in particular (\ref{br}). 
The following simple fact will be used various times.
\begin{remark}
\label{beta}
{\it Let $\nS:\YS\to X$ be the   normalization of $X$ at one node (i.e. $S=\{n\}$).
 Let $M\in \Pic \YS$ be such that $h^0(M)\geq 2$.
If $h^0(M(-q_1-q_2))=h^0(M)-2$,  every $L\in \FM$ satisfies
$h^0(X,L)=h^0(\YS,M)-1$.}

To prove it, pick $L\in \FM$ and consider the cohomology sequence
\begin{equation}
\label{C11}
0\to H^0(X,L) \stackrel{\alpha}{\la} H^0(\YS,M)\stackrel{\beta}{\la} k\to H^1(X,L) \to H^1(\YS,M)\to 0
\end{equation}
(associated to  (\ref{S1})).
It suffices to show that $\beta$ is  non zero.
The assumption  $h^0(M(-q_1-q_2))=h^0(M)-2$ implies that 
$h^0(M(-q_h))=h^0(M)-1$ for $h=1,2$;    hence 
$M$ has a section $s$ vanishing at $q_1$ but not at $q_2$;
but then $\beta (s)\neq 0$.
\end{remark}

\begin{nota}{}
\label{WMS}
Let  $S\subset\sing$ and consider the
 partial normalization $\YS\to X$. Fix
a finite set $S'$ of points of $X$   (usually $S'\subseteq S$). For any   $M\in \Pic \YS$
set
\begin{equation}
\label{WMSeq}
W_M(X,S'):=\{L\in \WM : \forall s\in H^0(X,L) \  \exists n\in S': \  s(n)=0 \}
\end{equation}
or equivalently (since $S'$ is finite)
\begin{equation}
\label{WMSeq2}
W_M(X,S'):=
\{L\in \WM : \exists n\in S': \  s(n)=0\  \   \forall s\in H^0(X,L)\}.
\end{equation}
If $S=\sing$ then
$W_M(X,S)$ is equal to the set of points in $\WM$ which  do not lie in  
$\amd (\dfib ^{\md})
$, where $\md =\mdeg M$. 
\end{nota}
\begin{lemma}
\label{h1}
Fix $\nS:\YS \to X$ and
let $M\in \Pic^{d}\YS$ be such that 
$h^0(\YS, M)= 1.$
\begin{enumerate}[(1)]
\item
\label{h1no}
If there exists $n_j\in S$ such that 
$h^0(\YS,M(-q^j_1))\neq h^0(\YS,M(-q^j_2))$, then $\WM =\emptyset$.
\item
\label{h1yes}
If $h^0(\YS,M(-q^j_1))=h^0(\YS,M(-q^j_2))$ 
for every $j$,  there are two cases.
\begin{enumerate}[(a)]
\item
\label{h1gen}
If $h^0(\YS,M(-q^j_h))=0$ for every $j$ and $h$, then  $\YS$ is connected and
 there exists a $L_M\in \FM$ such that
 $ \WM=\{L_M\}$ and $h^0(L_M)=1$. Moreover $\WMS =  \emptyset$
(hence $L_M\in \aXd(\dfib ^d) $).
\item
\label{h1all} If  
 there exists  $j$ for which $h^0(\YS,M(-q^j_1))=h^0(\YS,M(-q^j_2))=1$,
then $\WMS = \WM$.
Moreover:  if $h^0(\YS,M(-q^j_h))=1$ for every $j$ then $\WM=\FM $;
otherwise  $\WM=\{L_M\}$.
\end{enumerate}
\end{enumerate}
\end{lemma}
\begin{proof} Let  $s\in H^0(M)$ be a nonzero section. In case (\ref{h1no}) we are assuming that (up to
switching the branches over $n_j$)
 $s(q^j_1)=0$ while $s(q^j_2)\neq 0$, so obviously $s$ does not descend to a section of any $L\in \FM$.

For case (\ref{h1gen}) 
suppose, by contradiction, that  $\YS=\coprod_1^{\gamma}Z_i$  is not connected. Then
$h^0(Y,M)=\oplus h^0(Z_i,M_{Z_i})=1$ so that   there is only one component, say $Z_1$ such
that $h^0(Z_1,M_{Z_1})\neq 0$. Pick $q=q^j_h\in Z_2$, then (as $h^0(Z_2,M_{Z_2})=0$) every
$s\in H^0(M)$ vanishes at $q$ so that
$h^0(M(-q))=h^0(M)=1$ contradicting the hypothesis.
So $Y$ is connected.  Now  any nonzero $s\in H^0(Y,M)$ satisfies
 $s(q_h^j)\neq 0$ for
$\jv$ and $h=1,2$. Let $c_j:=s(q_2^j)/s(q_1^j)\in k^*$ and $\mc=(c_1,\ldots, c_{\dd}$);
then $\mc$ does not depend on the choice of $s$, as $h^0(M)=1$.  Using the construction 
of \ref{Lc}
set $L_M=L^{(\mc)}$;
we get $\WM =\{L_M\}$ 
and obviously $s$ descends to a section of $L_M$ that does not vanish at any $n_j$.
So, $W_M(X,S)$ is empty, and by construction, $h^0(X, L_M)=1$.

In case (\ref{h1all}), it is clear that for every $L\in \WM$ and $s\in H^0(L)$
 we have $s(n_j)=0$, hence $\WMS = \WM$. The last sentence is proved similarly.
\end{proof}

\begin{lemma}
\label{d1}
Let $\nS:\YS\to X$ be the   normalization of $X$ at one node (i.e. $S=\{n\}$).
Let $M\in \Pic^{d}\YS$ be such that 
$h^0(\YS, M)\geq  2$. 

Then $\WM =\FM$ and the following cases occur.

\begin{enumerate}[(1)]
\item
\label{d1gen} If $h^0(M(-q_1-q_2))=h^0(M)-2$ then $\WMS=\emptyset$ and
$h^0(L) = h^0(M)-1$ for every $L\in \FM$.
\item
\label{d1nobp}
If $h^0(M(-q_1-q_2))=h^0(M(-q_h))=h^0(M)-1$ for $h=1,2$ 
then  $\YS$ is connected and $\WMS=\WM \smallsetminus \{L_M\}$
for a uniquely determined $L_M\in \WM$
(hence $L_M\in \aXd(\dfib ^d) $). Moreover $h^0(L_M)=h^0(M)$ 
while  for every $L\in \WM - \{L_M\}$ we have $h^0(L)=h^0(M)-1$.
\item
\label{1bp}
If $h^0(M(-q_1))=h^0(M)-1$ and $h^0(M(-q_2))=h^0(M)$ then
$\FM=\WMS$. Moreover $h^0(L)=h^0(M)-1$ for every $L\in \FM$.
\item
\label{2bp}
If $h^0(M(-q_1))=h^0(M(-q_2))=h^0(M)$ then
$\FM=\WMS$. Moreover $h^0(L)=h^0(M)$ for every $L\in \FM$.
\end{enumerate}
\end{lemma}
\begin{proof} 
Pick $L\in \FM$ and consider the cohomology sequence (\ref{C11}).
%%$$
%%0\to H^0(X,L) \stackrel{\alpha}{\la} H^0(Y,M)\stackrel{\beta}{\la} k\la H^1(X,L) \la H^1(Y,M)\to 0
%%$$
It yields that $\alpha(H^0(X,L))$ has codimension at most $1$, i.e. that $h^0(L)\geq h^0(Y,M)-1\geq 1$
so that $\WM=\FM$. 
We shall omit the subscript $S$ during the proof.

In case (\ref{d1gen}), $H^0(Y,M(-q_1-q_2))$ has codimension $2$ hence
$\alpha(H^0(X,L))$ cannot be contained in it. Therefore $H^0(X,L)$
contains sections that do not vanish at $n$.
The rest has been proved in remark~\ref{beta}.

For the remaining cases, note that every section of $H^0(M(-q_1-q_2))$ descends to a section of every $L\in \FM$.

Case (\ref{d1nobp}).
To show that $Y$ is connected, suppose by contradiction that $Y=Y_1\coprod Y_2$, then (say)
$q_1\in Y_1$ and $q_2\in Y_2$ and $h^0(M)=h^0(Y_1,M_1)+h^0(Y_2,M_2)$
(denoting $M_i=M_{Y_i}$).  Furthermore 
$$
h^0(M_1)+h^0(M_2)-1=h^0(M)-1=h^0(M(-q_1))=h^0(M_1(-q_1))+h^0(M_2)
$$
hence $h^0(M_1(-q_1))=h^0(M_1)-1$. Similarly, $h^0(M_2(-q_2))=h^0(M_2)-1$.
But then 
$h^0(M(-q_1-q_2))=h^0(M_1(-q_1))+h^0(M_2(-q_2))=h^0(M)-2$
which is a contradiction.

Now, there exists $s\in H^0(M)$ such that $s(q_h)\neq 0$ for $h=1,2$.
Thus
\begin{equation}
\label{split}
H^0(M)=H^0(M(-q_1-q_2))\oplus k  s.
\end{equation}
Set $c=\frac{s(q_2)}{s(q_1)}$ and let $L_M=L^{(c)}$ (as in \ref{Lc}). The $s$ descends to a section
$\ov{s}\in H^0(L_M)$ such that $\ov{s}(n)\neq 0$. Hence $L_M\not\in \WMS$ and $h^0(L_M)=h^0(M)$.
Now, $L_M$ is uniquely determined, indeed if $s'\in H^0(M)$ is another section such that
$s'(q_h)\neq 0$ for $h=1,2$, then by (\ref{split}) $s'=a t + bs$ for $t\in H^0(M(-q_1-q_2))$ and $a,b\in k$
with $b\neq 0$. Thus $c=\frac{s'(q_2)}{s'(q_1)}$. This  proves that for every $L\in \WM$ such that
$L\not\in \WMS$ we have $L=L_M$.

In case (\ref{1bp}), $H^0(M(-q_1-q_2)=H^0(M(-q_1))$ and these are the only sections that can be pull backs of
sections of any $L\in \FM$.
Case (\ref{2bp}) is obvious.
\end{proof}

\begin{cor}
\label{00} %For every $X$ we have $W_{(0,\ldots,0)}(X)=\{\O_X\}$.
$W_{(0,\ldots,0)}(X)=\{\O_X\}$ for every  connected, nodal curve  $X$.
\end{cor}
%\begin{proof}
%Of course $\O_X\in \Wmd$. Let $\nu:Y\to X$ be the normalization of $X$.
%We can factor $\nu $ as follows
%$$\nu:Y\la Y_Y \stackrel{\nu_T}{\la} X$$
%with $Y_T$   connected of compact type.
%Thus $\Pic Y\cong \Pic Y_T$ and  $h^0(Y_T,\O_{Y_T})=1$ because $Y_T$ is
%connected. 
%Furthermore it is easy to see that $W_{(0,\ldots,0)}(Y_T)=\{\O_{Y_T}\}$ (using %\ref{h1} and \ref{d1}). To pass from $Y_T$ to $X$ we glue one pair of corresponding branches at the time; in each step
%we are in the situation of Lemma~\ref{h1} (\ref{h1gen}). Hence there exists a %unique $L\in \Pic^{(0,\ldots,0)}X$ such that $\nu_T^*L=\O_{Y_T}$,  $L=\O_X$ and we are done. 
%\end{proof}

\subsection{Divisors imposing independent conditions}
Let $\YS\to X$ be some partial normalization of $X$ and let $M\in \Pic \YS$. 
The goal of this subsection is to bound the dimension of the locus of $L\in W_M(X)$ 
which are not contained in the image of
the Abel map (i.e. with the notation of
\ref{WMS} the dimension of $\WMS$).
The easy cases,    $h^0(\YS,M)=1$ or  $\#S=1$, are dealt with by Lemmas~\ref{h1} and \ref{d1}.
To treat the general case we  introduce the following.

\begin{defi}
\label{ic}
Let $Y$ be a nodal curve (possibly   not connected).
Let $M\in \Pic Y$ and
let $E$ be a Cartier divisor on $Y$.
\begin{enumerate}[(A)]
\item
\label{ica} We say that $E$ is {\it admissibile for} $M$ if
for every subcurve $V\subseteq Y$ we have
$0\leq \deg_VE\leq h^0(V,M_V)$ (in particular,  $E$ is effective).
\item
\label{}
We say that {\it $E$ imposes independent conditions on $M$} if $E$ is admissible for $M$ and if 
$h^0(V, M(-E)_V)=h^0(V,M_V)-\deg_VE$
for every subcurve $V\subseteq Y$.
\item
\label{AMR} Let $R\subset Y\smallsetminus Y_{\text{sing}}$, we denote  by ${\mathcal A}(M,R)$
the set of all admissible divisors for
$M$ with support contained in $R$.
\end{enumerate}
\end{defi}
\begin{remark}
\label{fad}  If $R$ in part (\ref{AMR}) is finite, then
   the set ${\mathcal A}(M,R)$
 is also  finite.
\end{remark}

Let $C$ be a smooth irreducible curve, Definition~\ref{ic} gives back the classical one.
Fix a  finite subset   $R\subset C$; then every admissible divisor  $E$
such that $\Supp E\subset R$ imposes independent conditions on the general $L\in \Pic^dC$.
More generally

\begin{lemma}
\label{icst} Let  $\nu:Y\to X$  be the normalization of $X$ and $R\subset Y$   a finite
subset. Let $\md \in \BXs$   and
 $M\in \Pic^{\md}Y$   a general point.
Then  every   divisor  $E\in {\mathcal A}(M,R)$ 
 imposes independent conditions on $M$.
\end{lemma}
\begin{proof}
By Remark~\ref{fad}, it suffices
to prove that   a fixed $E$ imposes independent conditions on the general $M\in \Pic^{\md}Y$.

Set as usual $Y=\coprod_{i=1}^{\gamma} C_i$.
Given $M$ and $E$ as in the statement, denote $M_i:=M_{C_i}$, $E_i:=E_{C_i}$ and $e_i=\deg_{C_i}E$.
Now, for any line bundle $N$ on $Y$ and any subcurve $V\subset Y$ we have
$
H^0(V,N)=\oplus_{C_i\subset V} H^0(C_i, N_{C_i}).
$
Therefore it suffices to prove that
$ 
h^0(C_i, M(-E)_{C_i})=h^0(C_i,M_i)-e_i,
$  for every  $i=1\ldots,\gamma$.
Since $M$ is general in $\Pic^{\md}Y=\prod \Pic^{d_i}C_i$, every $M_i$ is general in $\Pic^{d_i}C_i$.
The fact that $\md$ is semistable implies that $d_i\geq p_a(C_i)-1\geq g_i-1$ ($g_i$ is the genus of
$C_i$) hence
$
h^0(C_i,M_i)=d_i-g_i+1.
$
Now  by (\ref{ica}) of \ref{ic} we have $e_i\leq d_i-g_i+1$ hence 
\begin{equation}
\label{degME}
\deg_{C_i}M(-E)=d_i-e_i\geq g_i-1.
\end{equation}
At this point, observe that $M_i(-E_i)$ is a general point in $\Pic^{d_i-e_i}C_i$ 
($E_i$ is fixed and $M_i$ is general in $\Pic^{d_i}C_i$) and hence
(by (\ref{degME}))
$$
h^0(C_i, M_i(-E_i))=d_i-e_i-g_i+1=h^0(C_i,M_i)-e_i
$$ as claimed.
\end{proof}
\begin{example}
\label{exad}
Let  $\nu:Y\to X$  the normalization of $X$  and
 $\md \in \BXs$.
 Then there exists a choice of branches $h:\{1,\ldots,\delta\}\to \{1,2\}$ 
such that the
divisor
$E=\sum_{j=1}^{\dd}q_{h(j)}^j$ is admissible for every  $M\in \Pic^{\md}Y$.
In fact, the construction of such an admissible divisor $E$ has  appeared in the proof of
\ref{GSb}. Recall that 
$\deg _{C_i}M(-E)=g_i-1$
for every $i=1,\ldots,\g$ (see (\ref{degM-E})).
\end{example}

For the next result we need some notation. Recall that
 $\nS:\YS \to X$ denotes the normalization  of $X$ at $S$.
Let $Z\subset X$ be a subcurve, we denote $Z_S:=\nu_S^{-1}(Z)$ 
the corresponding subcurve in
$\YS$, so that $Z_S$ is the  normalization of $Z$ at $S\cap Z_{{\text{sing}}}$.
Obviously   every subcurve of $\YS$ is of the form $Z_S$ for a unique $Z\subset X$.
We shall often simplify the notation  by setting
$
H^0(Z_S, M):=H^0(Z_S,M_{Z_S}).
$

\begin{prop}
\label{nvr}
Fix $\nS:\YS \to X$ as above.
Let $M\in \Pic\YS$ be such that
$h^0(\YS,M)\geq \dS$, and assume that for every $Z_S\subsetneq \YS$,
\begin{equation}
\label{ZS}
h^0(Z_S,M_{Z_S})\geq 1+\#(S\cap Z_{{\text{sing}}}).
\end{equation}
If every 
$E\in {\mathcal A}(M,\nS^{-1}(S))$
imposes independent conditions on $M$,
then 
\begin{displaymath}
\dim W_M(X,S)\leq\left\{ \begin{array}{ll}
\dS -\gS -1 &\text{ if } h^0(M)= \dS\\
\dS -\gS &\text{ if } h^0(M)\geq  \dS+1.\\
\end{array}\right.
\end{displaymath}
\end{prop}

\begin{proof}
We set
$l=h^0(\YS,M)$. 
By hypothesis, for every $q\in \nu^{-1}(S)$
\begin{equation}
\label{q1}
h^0(\YS, M(-q))=l-1,
\end{equation}
indeed by (\ref{ZS}) every such $q$  is admissible for $M$. 
Let $n\in S$ and set $\nu^{-1}(n)=\{q_1,q_2\}$.
Suppose $l=1$;
then $\dS=1$. 
By (\ref{q1}) applied to $q_1$ and $q_2$,
we are in case (\ref{h1gen}) of  Lemma~\ref{h1}. Hence 
 $\WMS=\emptyset$ and  we are done.

From now on, we assume $l\geq 2$.
Let $E=q_1+q_2$, then $E$ is admissible, i.e. 
$
\deg_{Z_S}E\leq h^0(Z_S, M_{Z_S})
$,
 for every subcurve $Z_S\subset \YS$. Indeed,  for every $Z_S$,
we have 
$
h^0(Z_S,M_{Z_S})\geq 1
$ by (\ref{ZS}). On the other hand
$
\deg_{Z_S}E\leq 2
$
and equality holds iff $Z_S$ contains both $q_1$ and $q_2$,
i.e. if and only if $Z$ is singular at $n$. In this case,  $h^0(Z_S,M_{Z_S})\geq 2$
by (\ref{ZS}). Therefore, by hypothesis,   for every $Z_S$
\begin{equation}
\label{q2}
h^0(Z_S, M(-q_1-q_2))=h^0(Z_S, M_{Z_S})-2.
\end{equation}

Assume $\dS=1$. 
By (\ref{q2}) we are in case (\ref{d1gen}) of lemma~\ref{d1}.
Thus $\WMS$ is empty and we are done.
We continue by induction on $\dS$. 

For every $\jv$ set
$S_j:=S\smallsetminus \{n_j\}$. For any  $\{j_1,j_2\}\subset \{1,\ldots, \dS\}$
\begin{equation}
\label{pair}
W_M(X,S)=\bigcup_{j=1}^{\dd}W_M(X,S_j)=W_M(X,S_{j_1})\cup W_M(X,S_{j_2})
\end{equation}
therefore it suffices to bound the dimension of   $W_M(X,S_j)$ for a chosen pair of  
values of $\jv$. 
We pick one of them
and simplify the notation by setting
$n=n_j$ and $T=S_j=S\smallsetminus \{n\}$.
We factor $\nS$
$$
\nS:\YS\stackrel{\nu_1}{\la}\YT \stackrel{\nT}{\la} X
$$
where  $\nT$ is the normalization of $X$ at $T$ and $\nu_1$ the normalization at the remaining node
$n$. We abuse notation by using the same names for points in $\YS$, $\YT$ and $X$ whenever the maps
are local isomorphims (e.g. $n$ denotes a node in $\YT$ and in $X$).
The following is the basic diagram to keep in mind
\begin{equation}
\label{diag}
\begin{array}{lcccr}
\Pic X & \stackrel{\nu_{T}^*}{\la}&\Pic \YT&\stackrel{\nu_1^*}{\la} &\Pic \YS \\
W_M(X,T)&\la&W_M(\YT)&\la&M\  \  \\
W_N(X,T)&\la &N &\mapsto&M\  \  
\end{array}
\end{equation}
where $N\in F_M(\YT)$;
since $l\geq 2$,
$
F_M(\YT)=\W_M(\YT).
$
By (\ref{q2}) and \ref{beta},
\begin{equation}
\label{h0N}
h^0(\YT, N)=l-1.
\end{equation}

\

\noindent
{\bf Case 1.} {\it The node  $n$  lies in two different irreducible components of $X$.}

By  Lemma~\ref{icind} part (\ref{extn})  (applied with $R=\nS^{-1}(S\smallsetminus n) $) every admissible divisor
$E_T$ on
$\YT$,  such that
$\supp E_T\subset \nT^{-1}(T)$, imposes independent conditions on $N$.
Therefore we can use induction ($\#T= \#S -1$) and obtain 
\begin{displaymath}
\dim W_N(X,T)\leq\left\{ \begin{array}{ll}
\dS -1-\gT -1 &\text{ if } h^0(\YT, N)= \dS-1\\
\dS -1 -\gT &\text{ if } h^0(\YT, N)\geq  \dS\\
\end{array}\right.
\end{displaymath}
i.e.  using (\ref{h0N})
\begin{displaymath}
\dim W_N(X,T)\leq\left\{ \begin{array}{ll}
\dS -\gT -2 &\text{ if } l-1= \dS-1\\
\dS -\gT -1&\text{ if } l-1 \geq  \dS .\\
\end{array}\right.
\end{displaymath}
If $n$ is not a separating node for $X$, then
$F_M(\YT)=\W_M(\YT)\cong k^*$ 
and $\gS=\gT$. Therefore, from diagram (\ref{diag}),
$
\dim W_M(X,T)\leq \dim W_N(X,T)+1
$.
So, using that $\dS -\gT -1=\dS-\gS$, we are done.

If  $n$ is separating, then $\gS = \gT +1$. On the other hand 
$\dim F_M(\YT)=0,$ 
hence $\dim W_M(X,T)\leq \dim W_N(X,T)$. Again, we are done.

\

\noindent
{\bf Case 2.} {\it The node $n$ 
 lies in only one irreducible component of $X$.}

Call $\ov{C}\subset X$ the component containing $n$ and 
 $C\subset \YS$ the component   containing both $q_1$ and $q_2$.
 We are in the
situation of Lemma~\ref{icind} 
part (\ref{intn}).
Therefore there exists a finite set $P\subset F_M(\YT)$ such that
for 
every $N\in \Pic \YT\smallsetminus P$, 
every admissible  $E$ supported on $\nT^{-1}(T)$  imposes independent conditions on $N$.
We can  use induction on 
 every $N\in W_M(\YT)$ such that $N\not\in P$. We obtain 
\begin{displaymath}
\dim W_N(X,T)\leq\left\{ \begin{array}{ll}
\dS -1-\gT -1 &\text{ if } h^0(\YT, N)= \dS-1\\
\dS -1 -\gT &\text{ if } h^0(\YT, N)\geq  \dS.\\
\end{array}\right.
\end{displaymath}
Consider diagram (\ref{diag}) and note that now  $\dim W_M(\YT)=\dim F_M(\YT)=1$.
Hence, away from the fibers over $P$, the dimension of every irreducible component of
$W_M(X,T)$ is at most 
\begin{displaymath}
\dim W_M(\YT)+\dim W_N(X,T)\leq\left\{ \begin{array}{ll}
1+\dd -\gamma -2 &\text{ if } l= \dd\\
1+\dd -\gamma -1&\text{ if } l\geq \dd+1\\
\end{array}\right.
\end{displaymath}
(using (\ref{h0N})) as wanted.

It  remains to bound the dimension of the fibers over every $N\in P$.
Now, set $n=n_1$ and
$T=\{n_2,\ldots, n_{\dS}\}$.

If $l\geq \dS +1$, i.e. if $h^0(\YT,N)\geq \dT+1$ then 
$$
\dim W_N(X)=\dim F_N(X)=\dT-\gT+1=\dS-\gS.
$$
The fiber of $W_M(X,T)\to W_M(\YT)$ over $N$ is obviously contained in $W_N(X)$,
hence it has dimension at most 
$\dS-\gS$ and 
 we are done.

Assume $\dS =l$.
If 
\begin{equation}
\label{GSNR}
h^0(\YT,N(-q^2_1-\ldots - q^{\dS}_1))=0,
\end{equation}
then, by \ref{GS} (\ref{GSW}), $W_{N}(X)$ has pure dimension equal to
$\dT -\gS=\dS-\gS -1$. Hence the dimension of the fiber
of $W_M(X,T)$ over $N$ is  at most $\dS-\gS -1$ and
we are done.

We shall complete the proof  by showing that  (\ref{GSNR}) holds for some choice of branches.
Assume $h^0(\YT,N(-q^2_1-\ldots - q^{\dS}_1))\geq 1$.

Observe that  $E:=\sum_{j=2}^{\dS}q_1^{j}+q^{\dS}_2$ is admissible for $M$.
Indeed, 
$
\deg_{Z_S}E\leq 1+ \#T\cap Z_{{\text{sing}}}
$ for every   $Z_S\subset\YS$;
hence, by (\ref{ZS}),
$$
\deg_{Z_S}E\leq 1+ \#T\cap Z_{{\text{sing}}}\leq 1 + \#S\cap Z_{{\text{sing}}}\leq h^0(Z_S,M).
$$
As $E$ is admissible, we have 
\begin{equation}
\label{absM}h^0(Y_S,M(-\sum_{j=2}^{\dS}q_1^{j}-q^{\dS}_2))=0,
\end{equation}
also, by Lemma~\ref{d1}, 
$$h^0(\YT,N(-q^2_1-\ldots -q^{{\dS}-1}_1- q^{\dS}_2))\leq 1 \  \text{ and } \  h^0(\YT,N(-q^2_1-\ldots -
q^{\dS}_1))= 1.$$ If
$
h^0(\YT,N(-q^2_1-\ldots -q^{{\dS}-1}_1- q^{\dS}_2))=1 
$
then, of course,
\begin{equation}
\label{abs}
h^0(N(-\sum_{j=2}^{\dS}q_1^{j}-q^{\dS}_2))= 1,
\end{equation}
which is impossible, by (\ref{absM}). Therefore 
$
h^0(\YT,N(-q^2_1-\ldots -q^{{\dS}-1}_1- q^{\dS}_2))=0 
$, i.e.  (\ref{GSNR}) holds for some choice of branches. The proof is complete.
\end{proof}
In the proof of Proposition~\ref{nvr} we used the following

\begin{lemma}
\label{icind}
 Let $\nu_1:\YS\la \YT$ be the partial normalization of $\YT$ at a unique node $n$.
Let $M\in \Pic\YS$ be such that for every subcurve $Z_S\subset \YS$
\begin{displaymath}
h^0(Z_S,M)\left\{ \begin{array}{ll}
\geq 2 &\text{ if } \nu_1^{-1}(n)\subset Z_S\\
\geq 1 &\text{ otherwise.}\\
\end{array}\right.
\end{displaymath}
Let $R$ be a finite set of smooth points of $\YS$.
Assume that  every   divisor  in  ${\mathcal A}(M,\nu_1^{-1}(n)\cup R)$
imposes independent conditions on $M$.
\begin{enumerate}[(i)]
\item
\label{extn}
If $n$ lies in two  irreducible components of $\YT$,  for any $N\in F_M(\YT)$, 
 every   divisor in ${\mathcal A}(N,\nu_1(R))$  imposes independent conditions on $N$.
\item
\label{intn} 
If $n$ lies in only one irreducible component of $\YT$,
 there exists a finite subset $P\subset F_M(\YT)$ 
such that for every $N\in F_M(\YT)\smallsetminus P$,
every divisor in ${\mathcal A}(N,\nu_1(R))$  imposes independent conditions on $N$.
\end{enumerate}
\end{lemma}
\begin{proof}
Let $\nu^{-1}(n)=\{q_1,q_2\}$.
Formula (\ref{q2}) holds  (with the same proof). 
For every  $Z_S\subset \YS$,  denote $Z_T:=\nu_1(Z_S)$. We have by (\ref{q2}) and \ref{beta}
\begin{equation}
\label{drop}
\text{if } \  \{q_1,q_2\}\subset Z_S \  \   \  \Longrightarrow \  \  h^0(Z_T, N_{Z_T})=h^0(Z_S,
M_{Z_S})-1
\end{equation}
and
\begin{equation}
\label{ZT=}
\text{if } \  \{q_1,q_2\}\not\subset Z_S \  \   \  \Longrightarrow \  \  h^0(Z_T, N_{Z_T})=h^0(Z_S,
M_{Z_S})
\end{equation}
because in this case $Z_S\cong Z_T$ via  $\nu_1$.
Thus for any $N\in F_M(\YT)$, the number $h^0(Z_T, N_{Z_T})$ 
depends only on
$M$, and not on the choice of $N$. Therefore the set 
 ${\mathcal A}(N,\nu_1(R))$ 
depends only on $M$.

Pick $E_T\in {\mathcal A}(N,\nu_1(R))$.
Denote $E_S:=\nu_1^*(E_T)$, and observe that $\nu_1$ is an isomorphism locally at every point in
$\supp E_S$.
Hence
\begin{equation}
\label{deg=}
\deg_{Z_S}E_S=\deg_{Z_T}E_T\leq h^0(Z_T,N)\leq h^0(Z_S,M).
\end{equation}
 Therefore $E_S$ imposes independent conditions on $M$,
i.e.
\begin{equation}
\label{icM}
h^0(Z_S,M(-E_S))=h^0(Z_S,M)-\deg_{Z_S}E_S.
\end{equation}
If $\{q_1,q_2\}\not\subset Z_S$,  $\nu_1$ induces an isomorphism $Z_S\cong Z_T$, hence by (\ref{deg=})
and (\ref{icM}) we get
$
h^0(Z_T,N(-E_T))= h^0(Z_S,M(-E_S))=h^0(Z_T,N)-\deg_{Z_T}E_T
$
as wanted. So we need only consider the case $\{q_1,q_2\}\subset Z_S$.

For part (\ref{extn}), let $q_1\in C_1$ and $q_2\in C_2$. Set
$e_i:=\deg_{C_i}E$ and $l_i:=h^0(C_i, M_{C_i})=h^0(C_i, N_{C_i})$.
Consider the usual  sequence
\begin{equation}
\label{C1T}
0\la H^0(Z_T,N(-E_T))\la H^0(Z_S, M(-E_S)) \stackrel{\beta}{\la} k \la \ldots
\end{equation}
If $E_T$ is such that 
$e_i\leq l_i-1$ for $i=1,2$ then $E_S+q_1+q_2$ imposes independent conditions on $M$. We get
$h^0(Z_S,M(-E_S-q_1-q_2))=h^0(Z_S,M(-E_S))-2$, hence 
$h^0(Z_T,N(-E_T))= h^0(Z_S, M(-E_S))-1$.
By (\ref{deg=}) and (\ref{icM}) we get
$$
h^0(Z_T,N(-E_T))= h^0(Z_S, M)-\deg_{Z_T}E_S-1=h^0(Z_T,N)-\deg_{Z_T}E_T
$$
as wanted.
Now, $E_T$ is admissible, hence $l_i\geq e_i$; so  only  two cases remain.

Case 1:  $e_1=l_1$ and $e_2=l_2-1$. Then $H^0(C_1,M(-E_S))=0$,  $h^0(C_2,M(-E_S))=1$
and $h^0(C_2,M(-E_S-q_2))=0$. Then all sections in $H^0(Z_S, M(-E_S))$ vanish at $q_1$
while there exist sections that do not vanish at $q_2$. Hence $\beta$ is surjective and
we are done.

Case 2: $l_i=e_i$ for $i=1,2$. Let $Z_T:=\nu_1(C_1\cup C_2)\subset Y_T$ .
By (\ref{drop}) 
$$
e_1+e_2=\deg _{Z_T}E_T\leq h^0(Z_T,N)=h^0(Z_S,M)-1\leq l_1+l_2-1
$$
which is possible only if at least one $e_i$ is less than $l_i$. So  Case 2 does not occur
and (\ref{extn}) is proved.

Now  part (\ref{intn}).
Call $C\subset \YS$ the component of $\YS$ containing both $q_1$ and $q_2$,
and $D:=\nu_1(C)$.
Set $e_D=\deg _{D}E_T=\deg _{C}E_S$;
and (by (\ref{drop}))
\begin{equation}
\label{}
l_D:=h^0(C,M)=h^0(D,N)+1
\end{equation}
so that $e_D\leq l_D-1$.
If $e_D\leq l_D-2$ then $E_S+q_1+q_2$ is admissible for $M$. Hence 
for every $Z_S\subset \YS$ we have $h^0(Z_S,M(-E_S-q_1-q_2))=h^0(Z_S,M(-E))-2$
so that (using remark \ref{beta})
\begin{equation}
\label{}
h^0(Z_T,N(-E_T))=h^0(Z_S,M(-E_S))-1=h^0(Z_T,N)-\deg_{Z_T}E_T.
\end{equation}
We are left with  case   $e_D=l_D-1$. Then $h^0(C, M(-E_S))=1$ and part (\ref{h1gen}) of
Lemma~\ref{h1} applies. We obtain that there exists a unique line bundle 
in $ \Pic  D$ which pulls back to $M(-E_S)_{C}$ and having $h^0=1$.
This in turn determines a  (unique) line bundle $N_D$ on $D$ 
which pulls back to $M_{C}$, and finally a unique line bundle on $\YT$ which pulls back to $M$ and restricts to 
$N_D$ on $D$. This last line bundle on $\YT$ is
 uniquely determined by 
$E_T$, so  we shall denote it by $\NE$.
Set 
$
P:=\{N\in F_M(\YT): N=\NE \text{ for some } E_T  \}
$.
We just showed  that
for any $N\in F_M(\YT)\smallsetminus P$,
 every $E_T\in {\mathcal A}(N,\nu_1(R))$  imposes independent conditions on $N$.
The finiteness of the set $P$
follows at once from the finiteness of  the set of  $E_T$'s.
\end{proof}

\begin{cor}
\label{WSgen}
Let $Y\to X$ be the normalization of $X$ and $S=\sing$.
If $\md \in \BX$ and $M\in \Pic^{\md}Y$ is a general point
then  $\dim \WMS \leq \dd - \gamma -1.$ 
\end{cor}
\begin{proof} If $M$ is general, $h^0(Y,M)=\dd$ by \ref{GSb}. Moreover, as $\md$ is stable, 
(\ref{ZS})  holds. Indeed
for every $Z\subset X$,   $Z^{\nu}=Z_S$
is the normalization of $Z$ and we have $d_Z\geq p_a(Z)= p_a(Z^{\nu})+\#Z_{\text{sing}}$;
hence $h^0(Z^{\nu}, M_{Z^{\nu}})\geq \#Z_{\text{sing}}+1$.
Finally, by Lemma~\ref{icst},  $M$ satisfies the assumption of Proposition~\ref{nvr}.

\end{proof}

\section{Irreducibility and dimension}
\subsection{Irreducible components}
We are ready to prove that $\Wmd$ is irreducible for every stable multidegree $\md$. 
This implies that, if $X$ is free from separating nodes, the theta divisor $\Theta(X)\subset \PXgb$ has one
irreducible component for every irreducible component of
$\PXgb$. If $X$ has some separating node this is  false (see \ref{ssred1} and 
\ref{ct}).
The stability assumption
on $\md$ is also essential, as one can
see from counterexample~\ref{ssred1}. 

If $|\md|\geq 1$ we shall use the Abel map $\amd$. If $|\md|\leq 0$, i.e. if
$g=0,1$ the Abel map is not defined so we need to treat this case separately, which will be done in the
following

\begin{lemma}
\label{g01} Let $X$ 
%% be a nodal connected curve of 
have genus $g\leq 0,1$;  let $\md\in \BX$.
Then 
\begin{displaymath}
\Wmd=\left\{ \begin{array}{ll}
\emptyset &\text{ if  } \md\neq (0,\ldots,0) \\
\{ \O_X\} &\text{ if  } \md = (0,\ldots,0) \  \  \text{ (hence   }$g=1$).\\
\end{array}\right.
\end{displaymath}
\end{lemma}
\begin{proof}
By hypothesis, $\forall \md \in \BX$ we have $|\md|=-1,0$ depending on whether $g=0,1$.
Recall that $X=\cup \ov{C_i}$ denotes the decomposition of $X$ in irreducible components.
Let $L\in \picX{\md}$ and suppose that
there exists a nonzero section $s\in H^0(X,L)$.  
Set
$$
Z^-:=\cup_{i:d_i<0}\ov{C_i};\  \  \  Z^0:=\cup_{i:d_i=0}\ov{C_i};\  \  \  
Z^+:=\cup_{i:d_i>0}\ov{C_i}.
$$
Note that $Z^-=\emptyset \Leftrightarrow \md =(0,\ldots,0)$. By contradiction, assume  
$Z^-\neq\emptyset$.
Then $s$ vanishes along a non empty subcurve $Z\subset X$ which contains $Z^-$.
Let $Z^C$ be the complementray curve of $Z$, so that $s$ does not vanish along any subcurve of $Z^C$.
Since for every $n\in Z\cap Z^C$ we have $s(n)=0$, the degree of $s$ restricted to $Z^C$ satisfies
\begin{equation}
\label{dzc}
d_{Z^C}\geq \#Z\cap Z^C.
\end{equation}
On the other hand,  $g\leq 1$ implies $p_a(Z^C)\leq 1$ hence the stability of $\md$ yields
$$
d_{Z^C}\leq p_a(Z^C)+\# Z\cap Z^C<\#Z\cap Z^C
$$
(cf. \ref{balrk}) a contradiction with (\ref{dzc}). Therefore $Z^-=\emptyset$; we obtain that, if $\Wmd\neq \emptyset$
then
$\md =(0,\ldots,0)$; in particular, $g=1$.
Now we conclude by Corollary~\ref{00}.
\end{proof}
Recall that for  $\md$ such that $\md \geq 0$ and $|\md|\geq 1$ we denote by $\Amd\subset
\picX{\md}$ the closure of the image of the Abel map $\amd$ (see \ref{abel}). If $\md\in \BX$ 
is such that $|\md|=-1,0$,  we denote $\Amd:=\Wmd$.

\begin{thm}
\label{irr} Let $X$ be a connected, nodal curve of arithmetic genus $g$.
Let $\md$ be a stable multidegree on $X$ such that $|\md|=g-1$. Then
\begin{enumerate}[(i)]
\item
\label{irr1}
$\Wmd = \Amd$, hence $\Wmd$ is irreducible of dimension $g-1$;
\item
\label{irr2}
the general $L\in \Wmd$ satisfies $h^0(X,L)=1$.
\end{enumerate}
\end{thm}
\begin{proof}
If $g=0,1$  the theorem follows from Lemma~\ref{g01}; so we  assume $g\geq 2$.
(\ref{irr2}) follows from (\ref{irr1}) and from \ref{GSb} (\ref{GSbh}).

Let $W$ be an irreducible component of $\Wmd$.
By \ref{pureW} we know that $\dim W= g-1$. 
We shall prove the Theorem by showing that if $ \Amd$ is not dense in $W$,
i.e. if $W\neq \overline{W\cap \im \amd}$,
then $\dim W \leq g-2$, and hence  $W$ must be empty.

Up to removing a proper closed subset of $W$, we can and will assume that 
$
W \cap  \im\amd = \emptyset .
$
Consider  $\nu:Y\to X$ the normalization of $X$, 
with $Y=\coprod_1^\g C_i$
and let $g_i$ be genus of $C_i$.
Recall that $g=\sumg+\dd -\gamma +1$.

We shall call $\rW$ the restriction to $W$ of the pull-back map $\nu^*$, so that
\begin{equation}
\label{rW}
\picX{\md}\supset W\stackrel{\rW}{\la}\rW (W)\subset \Pic^{\md}Y=\prod_{i=1}^\g \Pic ^{d_i}C_i.
\end{equation}
We shall bound the dimension of $W$ by analyzing $\rW$.

To say that  $L\in \picX{\md}$ does not lie in the image of
$\aXd:\dfib ^d\to \Pic X$ is to say that 
$L$ does not admit any section whose zero locus is contained in $\dfib$. In other words,
setting $S=\sing$, we have
$L\in \WMS$ (cf. \ref{WMS}). Therefore for every $M$  in $\rho(W)$ we have 
$$
\rho^{-1}(M)\subset \WMS\subset \WM.
$$
From now on,  $M$ is a general point in $\rho(W)$.
The proof is divided into four cases.

\noindent
{\bf Case I.} {\it $\dim \rho (W)\leq \sumg -2$.}

It suffices to use  that
$\dim\rho^{-1}(M)\leq \dim \FM =\dd-\gamma +1$. Then
$$
\dim W \leq \dim \rho(W)+\dim \FM\leq \sumg -2+\dd-\gamma +1=g-2.
$$

\noindent
{\bf Case II.} {\it $\dim \rho (W)= \sumg $.}

Now $\rho$ is dominant, so that $M$ is general in $\Pic^{\md}Y=\prod _1^{\gamma}\Pic^{d_i}C_i$.
Then we can apply Corollary~\ref{WSgen} which yields
$ \dim \WMS  \leq \dd -\gamma -1$, 
 and hence
$$
\dim W \leq \dim \rho(W)+\dim \WMS\leq \sumg +\dd -\gamma -1=g-2.
$$

\begin{remark}
\label{proj}
From now on we shall assume $\dim \rho (W)= \sumg -1$.

Denote $\pi_i:\prod
_{i=1}^{\gamma}\Pic^{d_i}C_i\to \Pic^{d_i}C_i$
the projection and $\rho_i:=\pi_i\circ \rho$ 
$$
\rho_i:W\la \rho(W)\la \rho_i(W)\subset \Pic^{d_i}C_i.
$$
As $\dim \prod _{i=1}^{\gamma}\Pic^{d_i}C_i=\sumg$ and $\dim \rho (W)= \sumg -1$, we get 
$$
\dim \rho_i(W)\geq g_i-1,\  \  \forall i
$$
and
there can be at   most one index $i$ for which   $\dim \rho_i(W)= g_i-1$.
\end{remark}

\noindent
{\bf Case III.} {\it $\dim \rho (W)= \sumg -1$ and $h^0(Y,M)\geq \dd+1$.}

We claim that we can apply \ref{nvr} to the general $M\in \rho(W)$. This would yield 
$ \dim \WMS  \leq \dd -\gamma$ 
so that we could conclude as follows:
$$
\dim W \leq \dim \rho(W)+\dim \WMS \leq g_Y -1+\dd -\g=g-2.
$$
To prove  that the hypotheses of \ref{nvr} hold,
observe  that (\ref{ZS}) follows from the fact that $\md$ is stable (see the proof of \ref{WSgen}).
To prove the remaining assumption 
 we argue by contradiction.
Assume that for some  admissible divisor $E$ with $\supp E\subset \nu ^{-1}(S)$ and 
$e:=\deg E$ we have
$$
h^0(Y,M(-E))\geq h^0(Y,M)-e+1
$$
for $M$ general in $\rho(W)$.
As $Y$ is the disjoint union of the $C_i$, we get 
$$
h^0(Y,M(-E))=\sum_{i=1}^{\gamma}h^0(C_i,M_i(-E_i))\geq \sum_{i=1}^{\gamma}(h^0(C_i,M_i)-e_i)+1
$$
where $E_i=E_{|C_i}$, $e_i:=\deg E_i$ and $M_i=M_{|C_i}$.
Therefore there exists at least one index, say $i=1$, such that
\begin{equation}
\label{1no}
h^0(C_1,M_1(-E))\geq h^0(C_1,M_1)-e_1+1.
\end{equation}
The fact that $E$ is admissible implies that $e_1\leq h^0(C_1,M_1)$. Now, as $d_1\geq g_1$,
there are two possiblities:

(a) $h^0(C_1,M_1)=d_1-g_1+1$;

(b) $h^0(C_1,M_1)\geq d_1-g_1+2$.

If (a) occurs,   $\rho_1:W\to \Pic^{d_1}C_1$ is dominant.  
In fact, by the assumption $h^0(M)\geq\dd +1$,  there exists an index $i\neq 1$
(say $i=2$)
such that $h^0(C_2,M_2)\geq d_2-g_2+2$, i.e. such that $M_2$ is a special line bundle on $C_2$.
Therefore   $\rho_2(W)$ cannot be dense in $\Pic^{d_2}C_2$.
By \ref{proj},  $\rho_1(W)$ is dense in $\Pic^{d_1}C_1$.
Therefore  we can apply Lemma~\ref{icst} (with $Y=X=C_1$ and $\md =d_1$),
getting that $E_1$ imposes independent conditions on $M_1$, a contradiction with (\ref{1no}).

%We can assume $e_1=d_1-g_1+1$ so that $h^0(C_1,M_1(-E))\geq 1$.Then the   map (denoting $E_1=E_{|C_1}$)
%\begin{equation}
%\label{}
%\begin{array}{lccr}
%u_{E_1}:& W^0_{d_1-e_1}(C_1)&\la &\Pic^{d_1}C_1 \\
%&N &\mapsto &N(+E_1)
%\end{array}
%\end{equation}
%is dominant. But this is impossible, for, by (\ref{dimW})
%$$\dim W^0_{d_1-e_1}(C_1)\leq \min\{d_1-e_1, g_1\}\leq \min\{g_1-1, g_1\}=g_1-1.$$

In case (b) we can assume $e_1=h^0(C_1,M_1)= d_1-g_1+2$. So $M_1$ is not a general point in $\Pic^{d_1}C_1$;
by \ref{proj},  $ \dim\rho_1(W)=g_1-1$.
Now (\ref{1no}) is $h^0(C_1,M_1(-E_1)\geq 1$. Consider the map
\begin{equation}
\label{uIII}
\begin{array}{lccr}
u_{E_1}:& W^0_{d_1-e_1}(C_1)&\la &\Pic^{d_1}C_1 \\
&N &\mapsto &N(+E_1).
\end{array}
\end{equation}
By what we said, $\im u_{E_1}$ dominates $\rho_1(W)$, hence
the variety $W^0_{d_1-e_1}(C_1)$ has
dimension at least $g_1-1$. This is impossible, since (by (\ref{dimW}))
$$
\dim W^0_{d_1-e_1}(C_1)\leq \min\{d_1-e_1, g_1\}\leq \min\{d_1-(d_1-g_1+2), g_1\}=g_1-2.
$$

\noindent
{\bf Case IV.} {\it $\dim \rho (W)= \sumg -1$ and $h^0(Y,M)= \dd$.}

If Proposition~\ref{nvr} applies,  we can argue  as for Case II and we are done.
Observe that
in order for \ref{nvr} to apply, it suffices to show that
 for every $i=1,\ldots,\gamma$, every divisor $E_i\in {\mathcal A}(M_i,\nu^{-1}(S)\cap C_i)$ 
%(cf. \ref{ic} (\ref{AMR}))  
imposes independent conditions on $M_i$.
Indeed this implies
that every $E\in {\mathcal A}(M,\nu^{-1}(S))$ imposes independent
conditions on $M$.
By \ref{proj} there are two possibilities.

\noindent
(a) $\rho_i(W)$ is dense in $\Pic^{d_i}C_i$ for every $i$.

\noindent
(b) There exists a unique index, say $i=1$, such that $\dim\rho_1(W)=g_1-1$, whereas for $i\geq 2$,
$\rho_i(W)$ is dense.

In case (a), $M_i$ is general in $\Pic^{d_i}C_i$, hence by \ref{icst} and by what we observed above we 
can use \ref{nvr} and we are done.

Now  case (b). We may
 assume that \ref{nvr} cannot be applied.
Let 
$
E:=\sum_{j=1}^{\dd}q_{h(j)}^j
$
be an admissible divisor for $M$ of the same type constructed in \ref{exad}
(with the same notation).
Recall  from \ref{exad} that $\deg_{C_i}M(-E)=g_i-1$, for all $i$.

If $E$ imposes independent conditions, i.e. if $h^0(Y,M(-E))=h^0(M)-\dd=0$,
   we can apply \ref{GS} (\ref{GSW}) and obtain that
$
\dim \WM = \g-\dd.
$
This is enough to  conclude:
\begin{equation}
\label{last}
\dim W \leq \dim \rho(W)+\dim\WM= \sumg -1+\dd -\gamma =g-2.
\end{equation}
So, assume that  $h^0(Y,M(-E))\geq 1$ .
We have that  $h^0(C_i,M_i(-E_i))=0$ if $i\geq 2$, whereas
 $h^0(C_1,M_1(-E_1))\geq 1$. As we said 
$\deg M_1(-E_1)=g_1-1$;  
we claim that
\begin{equation}
\label{h0M-E}
h^0(C_1,M_1(-E_1))= 1.
\end{equation}
To prove it we argue as for Case III (b). 
Consider the map analogous to (\ref{uIII}):
$$
u_{E_1}^1:  W^1_{g_1-1}(C_1) \la  \Pic^{d_1}C_1 
$$
mapping $N$ to $N(E_1)$.
Now   $\dim W^1_{g_1-1}(C_1)\leq g_1-3$ (well known);
therefore, $u_{E_1}^1$ cannot dominate $\rho_1(W)$, 
whose dimension is $g_1-1$.  So (\ref{h0M-E}) is proved.

It is trivial to check that we can assume, for a suitable $q\in \Supp E_1$,
that  $E_1=E_1'+q$ with $E_1'$ imposing independent conditions on $M_1$, i.e.
$$
h^0(C_1,M_1(-E'_1))=1
$$
so that $q$ is a base point of $M_1(-E'_1)$.
Therefore
$$
h^0(Y,M(-E))= 1
$$
and there exists a point $q\in E_1$ such that, setting $E'=E-q_1$,
the divisor $E'$ imposes independent conditions on $M$.
Now let $n$ be the node of $X$ of which the point $q_1$ is a branch,  
let $S'=S\smallsetminus n$; thus
$E'$ is supported on $\nu^{-1}(S')$.
Let
$\nu_n: X'\to X$ be the normalization of $X$ at $n$, so that we can factor $\nu$
$$
Y\stackrel{\nu'}{\la}X'\stackrel{\nu _n}{\la} X
$$
and $\nu '$ is the normalization of $X'$. Of course, $X'$ has $\delta '=\dd -1$
nodes and $h^0(Y,M)=\dd '+1$.
As $E'$ imposes independent conditions on $M$, we can apply \ref{GS} with respect to
$\nu':Y\la X'$. This gives us that $ W_M(X')=F_M(X')$
and, for a general $L'\in W_M(X')$
\begin{equation}
\label{L'}
h^0(X',L')=h^0(Y,M)-\dd'=1.
\end{equation}
Consider the following diagram
\begin{equation}
\label{diag1}
\begin{array}{lcccr}
\Pic X & \stackrel{\nu_n^*}{\la}&\Pic X'&\stackrel{(\nu ')^*}{\la} &\Pic  Y \\
W_M(X)&\la&F_M(X')&\la&M.\  \  
\end{array}
\end{equation}
Observe that   $n$ is not a separating node  of $X$ (otherwise, by \ref{empty},
$\BX$ is empty and there is nothing to prove).
Hence $\nu_n^*$ is a $k^*$-fibration and
$$
\dim F_M(X')=\dd'-\gamma +1=\dd -\gamma.
$$
We now claim that the fiber $W_{L'}(X)$  of $W_M(X)$ over the general point $L'\in F_M(X')$
 has dimension $\leq 0$.
By (\ref{L'}) we are in the situation of Lemma~\ref{h1}, which tells us that the only case when
$\dim W_{L'}(X) =1$ is when $L'$ has a base point in each of two branches of $n$.
Now this does not happen. Indeed, if $i\geq 2$  $M_i$ is general and hence has no base point over $\sing$;
on the other hand $M_1$ varies in a codimension $1$ subset of $\Pic ^{d_1}C_1$
hence it has at most one base point over $\sing$; therefore we can apply
Lemma ~\ref{nbp}. 

Concluding:
$\dim W_M(X)\leq \dd -\gamma$.
Arguing as in (\ref{last}) we are done.
 \end{proof}

\begin{example}
\label{ssred1} The Theorem fails if we assume $\md$ semistable.
The simplest instance of $\md\in \BXs$ with $\Wmd$  reducible  is that of a curve of compact type,
$X=C_1\cup C_2$, where $C_i$ is smooth of genus $g_i$,  $\#C_1\cap C_2=1$ and  $\md=(g_1-1,g_2)$ 
(note that $\md$ is strictly semistable by \ref{empty}).
Then
$$
\Wmd=(W_{g_1-1}(C_1)\times \Pic^{g_2}C_2) \cup (\Pic^{g_1-1}C_1\times \Theta_{q_2}(C_2))
$$
where $q_2\in C_2$ is the point  over the node and 
$\Theta_{q_2}(C_2):=\{L\in \Pic^{g_2}C_2
:h^0(C_2,L(-q_2))\ne 0\}.$ 
The interested  reader will easily construct similar, more interesting, examples on  curves 
 not
of compact type. 
\end{example}

\subsection{Dimension of the image of the Abel map}

\begin{prop}
\label{dimA} Let $X$ be a curve of genus $g\geq 2$.
Let $\md\in \Z^{\g}$ be a non-negative multidegree such that $|\md|=g-1$.
If $\md$ is semistable, then 
\begin{enumerate}[(a)]
\item
\label{hb}
the general $L\in \Amd$ satisfies $h^0(X,L)=1$;
\item
\label{dimb}
$\dim \Amd =g-1$.
\end{enumerate}
Conversely, if $\md$ is not semistable, then
\begin{enumerate}[(A)]
\item
\label{hnb}
for every $L\in \Amd$ we have $h^0(X,L)\geq 2$;
\item
\label{dimnb}
$\dim \Amd \leq g-2$.
\end{enumerate}
\end{prop}
\begin{proof}
If $\md$ is stable,  by Theorem~\ref{irr} we know that $\Amd =\Wmd$, 
 $\dim \Amd =g-1$ (by \ref{pureW}) and  that the general point $L\in
\Amd$ has
$h^0(X,L)=1$.
So, for the first half of the statement, we need to consider the  case  where $X$ is reducible and 
 $\md$ semistable but not stable.
Thus, there exists
a decomposition
 $X=V\cup Z$, where $V$ and $Z$ are subcurves of respective arithmetic genus $g_V$ and $g_Z$,
 such that $V$ is connected,
\begin{equation}
\label{bBI}
d_V= g_V-1 \  \text{ and } \  d_Z= g_Z+\dS -1,
\end{equation}
where  $S:=V\cap Z$ and $\dS:=\#S$. 

Observe that, since $\md \geq 0$, we get  
$
g_V\geq 1 $.
By (\ref{genusS}) we have
\begin{equation}
\label{genus2}
g=g_V+g_Z+\dS-1.
\end{equation}
Let $L$ be a general point in $\Amd$;  we can assume that $L$ is a line bundle on $X$
of type $L=\O_X(D)$ where $D$ is an effective divisor of multidegree $\md$ supported on the smooth locus of $X$.
Consider the restrictions $L_V$ and $L_Z$ of $L$ to $V$ and $Z$; we have
$h^0(V,L_V)\geq 1$. On the other hand
$h^0(Z,L_Z)\geq d_Z-g_Z+1=\dS$ (by Riemann-Roch and (\ref{bBI})); moreover equality holds for a general 
$L_Z\in \Pic^{\md_Z}Z$,
by Corollary~\ref{GSc}. 
Denote the partial normalization of $X$ at $S$ by
$$
\nS: \YS=V\coprod Z \la X
$$
and note that $\Pic^{\md}\YS=\Pic^{\md_V}V\times \Pic^{\md_Z}Z$.
Set $M=\nS^*L=(L_V,L_Z)$.
Then   for $L$ general 
%(i.e. for $L_Z$ general in $\Pic^{\md_Z}Z$)
\begin{equation}
\label{hM1}
h^0(\YS, M)=h^0(V,L_V)+h^0(Z,L_Z)= \dS +1
\end{equation}
hence by Proposition~\ref{GS} (\ref{GSh}), which we can apply by Lemma~\ref{GSb}(\ref{indb}),
we obtain that $h^0(X,L)=h^0(\YS, M)-\dS=1$.

Now we compute $\dim \Amd$ using
induction on the number of irreducible components of $X$. 
The case $X$  irreducible has already been settled.
Assume $X$ reducible; by what we said above, the pull back map $\nS^*$ restricted to $\Amd$
gives a dominant rational map (denoted by $\rho$)
$$
\Amd \stackrel{\rho}{\la } A_{\md _V}(V)\times \Pic^{\md_Z}Z.
$$
Now recall that 
$|\md_V |=g_V-1\geq 0$ by (\ref{bBI})  and $\md_V $ is semistable on $V$ because $\md$ is semistable on
$X$ (cf. \ref{balrk}).
Furthermore $V$ has fewer components than $X$, 
hence we can use
induction to conclude that 
$\dim A_{\md _V}(V)= d_V=g_V-1$.  
 
If $M$ is a general point in the image of the above map $\rho$, then by (\ref{hM1}) and \ref{GS} (\ref{GSW}),
we see that $\WM = \FM$. We claim that $\WM\subset \Amd$. Indeed recall that
$M=\nS^*\O_X(D)$ 
with $\Supp D\subset \dfib$, hence there exists an $L\in \WM$ (namely, $L= \O_X(D)$)
admitting a section that
does not vanish at any node of $X$. Therefore the same holds for
every line bundle in a dense open subset of
 $\WM$ (which is irreducible, being equal to $\FM$). This shows that $\WM \subset \Amd$.
Therefore $\rho ^{-1}(M) =\FM$ and
$$
\dim \Amd = g_V-1+g_Z+\dS - \gS +1=g-1.
$$

Conversely assume that $\md$ is not semistable. Then there exists a decomposition
$X=V\cup Z$, where (as before) $V$ and $Z$ are subcurves of  genus $g_V$ and $g_Z$
 such that 
\begin{equation}
\label{noBI}
d_V\leq g_V-2 \  \text{ and } \  d_Z\geq g_Z+\dS
\end{equation}
where  $S:=V\cap Z$ and $\dS:=\#S$.
Notice that  $g_V\geq 2$ (as $\md \geq 0$).

We use  the same notation as before.
Let $L$ be a general point in $\Amd$, so  $L$ is 
of type $L=\O_X(D)$ with $D\geq 0$  supported on 
$\dfib$. We have
$h^0(V,L_V)\geq 1$ and $h^0(Z,L_Z)\geq d_Z-g_Z+1\geq \dS +1$.

Consider
$
\nS: \YS=V\coprod Z \la X
$
and set $M=\nS^*L=(L_V,L_Z)$.
We have 
\begin{equation}
\label{hM}
h^0(\YS, M)=h^0(V,L_V)+h^0(Z,L_Z)\geq \dS +2
\end{equation}
hence by \ref{GS} (\ref{eff})
$
h^0(X,L)\geq 2,
$
proving part (\ref{hnb}).
To compute $\dim \Amd$ 
consider again the  rational map
$$
\Amd \stackrel{\rho}{\la } A_{\md _V}(V)\times \Pic^{\md_Z}Z.
$$
Since $\dim A_{\md _V}(V)\leq d_V$ (by Lemma~\ref{Ad}) we get
 $$
\dim \Amd \leq d_V+g_Z+\dim \WM\leq g_V-2+g_Z+\dim \WM
$$
using (\ref{noBI}) for the last inequality.
Thus
$$
\dim \Amd \leq g_V-2+g_Z+\dS - \gS +1=g-2.
$$
This proves (\ref{dimnb}) and we are done.
\end{proof}
From the proof, it is clear  that the farther is $\md$ from  semistable, the smaller is the dimension of $\Amd$. The following fact will be useful later on.

\begin{cor}
\label{bpf} Let $R\subset X$ be a finite set of nonsingular points of $X$
and $\md\in \BXs$. Then the general $L\in \Amd$ has no base point in $R$.
\end{cor}
\begin{proof}
It obviously suffices to assume $\#R=1$, so let $R=\{q\}$.
If $L$ is general in $\Amd$ we can assume that $L\in\im \amd$. Set $L'=L(-q)$ and
$\md ':=\mdeg L'$. If $q$ is a base point of $L$, then $L'\in \im\alpha^{\md'}_X$.
Therefore, if the general $L\in \Amd$ has a base point in $q$,
the map
\begin{equation}
%\label{}
\begin{array}{lccr}
 & \im\alpha^{\md'}_X&\la &\Amd \\
&L' &\mapsto &L'(q)
\end{array}
\end{equation}
must be dominant. But this is not possible, as $\dim \Amd =g-1$ by \ref{dimA}, whereas
obviously
$\dim \im\alpha^{\md'}_X\leq |\md'|=g-2$.
\end{proof}

\section{The compactified theta divisor}
\label{comp}
\begin{nota}{}
\label{blow}
Let $X$ be a connected nodal curve, $S\subset \sing$,  
 $\dS:=\#S$ and $\nS:Y_S\to X$   the normalization of $X$ at  
$S$.
Let 
\begin{equation}
\label{XSh}
\hXS=\YS\cup _{i=1}^{\dS}E_i
\end{equation}
be the  connected, nodal curve obtained by ``blowing up" $X$ at $S$, so that $E_i\cong \pr{1}, \forall i$
and $E_i$ is called an exceptional component of $\hXS\to X$
(where this map is the contraction of all the exceptional components of $\hXS$).
%Denote $\sigma :\hXS\to X$ the map contracting all the $E_i$. 
We shall usually denote by   $\hM$  a line bundle on $\hXS$ and by $M\in 
\Pic\YS$ its restriction to $\YS$.
\end{nota}
\subsection{The compactified  Picard variety}
\begin{nota}{}
\label{descgen} In what follows we shall recall what the points of $\PXgb$ parametrize,
and give a stratified description of it (in \ref{Pstr}); our notation is  that of \cite{cner}.
 There is more than one
place where  details and proofs can be found, even though some terminology may be different from ours.
We  refer to   \cite{alex} for a unifying account and other references.
 
To begin with, using the notation of \ref{blow},   
 the compactified  Picard variety, or compactified jacobian, $\PXgb$, in degree $g-1$,
 parametrizes  equivalence classes of stable line bundles of degree $g-1$ on  the curves $\hXS$ 
  as $S$ varies among all subsets of $\sing$.

Let us define stable line bundles 
and   the   equivalence relation among
them.
For every $S\subset \sing$ consider the blow up of $X$ at $S$,
$\hXS=\YS\cup _{i=1}^{\dS}E_i$ (cf.
(\ref{XSh})). A stable line bundle 
$\hM\in
\Pic^d\hXS$ is such that, setting $M:=\hM_{Y_S}$,   properties (1) and (2) below hold:

\

\noindent (1) $\mdeg M \in \BYS$;

\noindent (2) $\deg_{E_i} \hM = 1$ for    $i=1,\ldots \dS$.

\

\noindent
We call  $\hM\in
\Pic^d\hXS$ {\it semistable} if it satisfies (2) as well as (1$'$)  below

\

\noindent (1$'$) \  $\mdeg M \in \BYSs$.

\

In other words, a line bundle on $\hXS$  is semistable (resp. stable) if its restriction to the complement of all the
exceptional components of $\hXS\to X$ has semistable (resp. stable)
multidegree.
 Two stable line bundles $\hM$
and
$\hM '$ on
$\hXS$ are defined to be equivalent  iff their restrictions,    $M$ and $M'$,
 to $\YS$ coincide.
\end{nota}

\begin{nota}{}
\label{notpt}
Thus, the points in $\PXgb$ are 
 in one-to-one correspondence with equivalence classes of
stable line bundles. Any such class is uniquely determined by $S$ and by $M\in \Pic \YS$
(provided that $\BYS$ is not empty),
therefore  points of $\PXgb$ will be denoted by pairs
$
[M,S]
$, where $S\subset \sing$ and $M\in \Pic^{\md}\YS$ with $\md \in \BYS$.
\end{nota}

\begin{nota}{}
\label{GIT} Although $\PXgb$ is constructed as a GIT-quotient, our terminology
``stable/ semistable line bundles" does not precisely reflect the GIT  stability/ semistability.
More precisely, 
denote $q_X:H_X\to \PXgb$ the GIT quotient defining $\PXgb$ (so that $H_X$ is a closed
subset in the GIT-semistable locus of some Hilbert scheme).
Note that $H_X$ contains  strictly GIT-semistable points,
unless $X$ is irreducible.
Our stable line bundles correspond to GIT-semistable points in $H_X$  having
 closed orbit.
\end{nota}

\begin{nota}{}
\label{S(d)}
For technical reasons we need to consider semistable multidegrees that are not stable. 
Let $\md\in \BYSs$ be a semistable multidegree of $\YS$;
a  node $n$  of $\YS$ is called  {\it destabilizing} for $\md$ if there exists a  connected
subcurve $Z\subset \YS$ such that $n\in Z\cap Z^C$ and $d_Z=p_a(Z)-1$
($Z^C=\ov{Y\smallsetminus Z}$).
We set 
\begin{equation}
\label{Sd}
\Sd:=\{n\in (\YS)_{\text{sing}}:\  n \text{ is destabilizing for } \md\}.
\end{equation}
Observe that
\begin{equation}
\label{}
\Sd=\emptyset \  \Leftrightarrow \  \md\in \BYS.
\end{equation}
We call $\Yd$ the normalization of $\YS$ at $\Sd$, so that we have
\begin{equation}
\label{Yd}
\Yd=Y_{S\cup \Sd}\stackrel{\nd}{\la}\YS\stackrel{\nS}{\la}X
\end{equation}
where $\nd$ is the normalization map.

Assume that $\md$ is strictly semistable, i.e. $\Sd$ is not empty. Then the dual graph of $\YS$ has
an orientation such that for every subcurve $Z\subset \YS$ such that $d_Z=p_a(Z)-1$,
all the edges between $\Gamma_Z$ and $\Gamma_{Z^C}$ go from $\Gamma_Z$ to $\Gamma_{Z^C}$
(by \ref{beau}).
Therefore, if we consider $\Yd$ 
and use the convention of \ref{he}, 
for every destabilizing node $n\in Z\cap Z^C$, we have 
$q_1^n\in Z$ and $q_2^n\in
Z^C$ (abusing notation by denoting  $Z=\nd^{-1}(Z)$ and $Z^C=\nd^{-1}(Z^C)$). We now introduce a divisor on $\Yd$:
\begin{equation}
\label{Td}
T(\md ):=\sum_{n\in \Sd}q_2^n \  \text{ and }\  \  \mt(\md):=\mdeg T(\md ).
\end{equation}
By construction, $\md-\mt(\md)$ is a stable multidegree for $\Yd$. Set
\begin{equation}
\label{dst}
\dst:=\md - \mt(\md) \in \Sigma (\Yd ).
\end{equation}
\end{nota}
The following statement  summarizes various  known facts about
$\pXgb$. The only novelty is that we
use line bundles on the partial normalizations of $X$, rather
than torsion free sheaves on $X$ 
(as in \cite{AK}, \cite{OS}, \cite{simpson}) or line bundles on the blow-ups of $X$ (as in
\cite{caporaso}). 
\begin{fact}
\label{Pstr}
$\pXgb$ is a connected, reduced, 
projective scheme of pure dimension $g$.
It  has a stratification
$$
\pXgb =\coprod_{\stackrel{\emptyset\subseteq S \subseteq\sing}{\md \in \BYS}}\PSd
$$
such that the following properties hold.
\begin{enumerate}[(i)]
\item
\label{eps}
For every $S\subset\sing$ and every  $\md\in \BYS$ there is a canonical isomorphism
(notation in \ref{notpt})
\begin{equation}
%\label{}
\begin{array}{lcr}
\Pic^{\md}\YS&\stackrel{\eSd}{\la} & 
\PSd\\
\  \  M &\mapsto &[M,S].\\
\end{array}
\end{equation}
%Note: we shall often use $\eSd$ to identify  $\PSd$ with $\Pic^{\md}\YS$.
In particular, if $\PSd\neq \emptyset$, then $\dim\PSd=g-\dS+\gS-1$.
\item
\label{} (More generally) For every $S\subset\sing$ and every  $\md\in \BYSs$ there is a canonical
surjective morphism $\eSd:\Pic^{\md}\YS\to P^{\dst}_{\Sd}$
(notation in \ref{S(d)}) which factors as follows
\begin{equation}
%\label{}
\begin{array}{lccccr}
\eSd: & \Pic^{\md}\YS&\stackrel{\tau}{\la} & 
\Pic^{\dst}\Yd&\stackrel{\epsilon^{\dst}_{\Sd}}{\la}&P^{\dst}_{\Sd}\\
\;&L &\mapsto &\nd^*L\otimes\O_{\Yd}(-\sum _{n\in \Sd}q_2^n)& &\\
\end{array}
\end{equation}
where $\tau$ is surjective with 
 fibers   $(k^*)^b$,
$b=\dd_{\Sd}-\g_{\Sd}+1$, and  $\epsilon^{\dst}_{\Sd}$ is
an isomorphism.
\item
\label{stcl} If $P^{\md '}_{S'}\subset \overline{\PSd}$ then  $S\subset S'$ and $\md \geq \md'$
(i.e. $d_i\geq d_i', \forall i=1,\ldots, \gamma$).

In such a case,
$\#\bigl((S'\smallsetminus S)\cap \ov{C_i}\Bigr) =d_i-d_i'$ (recall that $X=\cup_1^{\g} \ov{C_i}$).

\item
\label{smooth}
Denote $\PXg$ the smooth locus of $\PXgb$, then
$$
\PXg=\coprod_{\md \in \Sigma(\widetilde{X})}P^{\md}_{\sep}\cong\coprod_{\md \in
\Sigma(\widetilde{X})}\Pic^{\md}\widetilde{X}
$$
where $\widetilde{X}\to X$ 
is the normalization at the separating nodes (cf. \ref{Ssplit})
and the isomorphism is the canonical one described in part (\ref{eps}).
\end{enumerate}
\end{fact}
Given the normalization of $X$ at all of its separating
nodes, $\widetilde{X}\to X$, recall from \ref{Ssplit} that $\widetilde{X}=\coprod_{i=1}^c
{X}_i$ denotes the connected components decomposition of
$\widetilde{X}$.
\begin{cor}
\label{irrP}
$\PXgb$ is irreducible if and only if 
for every $i=1,\ldots, c$ either ${X}_i$ is irreducible,
or
every irreducible component $C$ of ${X}_i$ meets  $\ov{{X}_i\smallsetminus C}$ in exactly $2$ points.
\end{cor}
\begin{proof}
 Assume first   $X=\widetilde{X}$. 
Then $\PXgb$ is irreducible if and only if $\#\BX=1$.
If $X$ is irreducible, then  $\BX=\{g-1\}$ so $\PXgb$ is irreducible.
If  every irreducible component $\ov{C}_i$ of $X$ meets its
complementary curve in $2$ points, calling $\ov{g}_i$
the arithmetic genus of $\ov{C}_i$,  we have
$ g-1=
%\sum_{i=1}^{\gamma} \ov{g}_i+\gamma -\gamma =
\sum_{i=1}^{\gamma}
\ov{g}_i.
$ Therefore
$\BX= \{ (\ov{g}_1,\ldots,
\ov{g}_{\gamma}) $\}, hence
$\PXgb$ is irreducible.

Conversely assume that
$X$ is reducible and has an irreducible
component, $\ov{C}_i$, such that 
$\dd_i:=\#\ov{X\smallsetminus \ov{C}_i}\geq 3$.
Then $X$ may be obtained as the special fiber of a family of
nodal curves $X_t$ having exactly two  irreducible components
intersecting  in $\dd_i$ points. 
Then  $\#\Sigma (X_t)=\dd _i-1\geq 2$ (cf. \ref{theta2}), hence
 $\ov{P_{X_t}^{g-1}}$  has at least $2$
irreducible components. Since  $\ov{P_{X_t}^{g-1}}$ specializes to
$\PXgb$ we get that $\PXgb$ has at least $2$ irreducible
components.  So, if $X$ has no separating node we are done.

In general, denote $\widetilde{b}:=\#\Sigma(\widetilde{X})$.  Then
 $\PXgb$ is irreducible if and only if  $\widetilde{b}=1$;
 by \ref{Ssplit} this is equivalent to  $\#\Sigma(X_i)=1$
 for every $i=1,\ldots,c$. Then the result follows by applying to each $X_i$ the first part.
\end{proof}

\begin{remark}
In combinatorial terms, consider the
graph  $\widetilde{\Gamma} _X$  obtained from $\Gamma_X$ by removing every loop and every
 separating edge. Then
$\PXgb$ is irreducible if and only if every vertex of $\widetilde{\Gamma}_X$
has valency (or degree) equal to either $0$ or $2$.
\end{remark}

\subsection{Stratifying the theta divisor.}
We shall now define the theta divisor of $\PXgb$ using
the stratification given above. A natural thing to do  is to
consider the irreducible strata, $\PSd$, of dimension $g$ of $\PXgb$, consider $\Wmd$ in such strata
and then take their closure.
Recalling Lemma~\ref{empty},  the $g$-dimensional strata are easily listed. 
First, denote
$\sep\subset \sing$ the set of separating nodes of $X$ and let
$
\widetilde{X} \to X
$
be the normalization of $X$ at $\sep$ (as in \ref{Pstr} (\ref{smooth})).
Thus $\widetilde{X}$ is a nodal curve having $c=\# \sep+1$ connected components.
Finally, set $\widetilde{b}=\#\Sigma(\widetilde{X})$. 
We have
\begin{lemmadef}
\label{thetadef} Let $X$ be a connected nodal curve. 
Using $\eSd$  of \ref{Pstr} (\ref{eps}) as an identification,
we define the theta divisor $\Theta(X)$ of $\PXgb$ as
$$
\Theta(X):=\overline{\bigcup_{\md \in \Sigma(\widetilde{X})}
W_{\md}(\widetilde{X})}\subset \PXgb.
$$
$\Theta (X)$ has $c\widetilde{b}$ irreducible components, all of dimension
$g-1$.
\end{lemmadef}
\begin{proof} 
If $X$ is free from separating nodes (i.e. $c=1$) the statement
 follows trivially from Theorem~\ref{irr}.
Otherwise, let  $\widetilde{X}=X_1\coprod\ldots \coprod X_{c}$
be the decomposition into connected components. Then $g=\sum_1^{c} p_a(X_i)$
and 
$$
W_{\md}(\widetilde{X})=\bigcup_{i=1}^{c}\bigr( W_{\md_i}(X_i)\times 
\prod_{\stackrel{j\neq i}{j=1,\ldots,c}}\Pic^{\md_j}X_j\bigl)
$$
where $\md_i$ denotes the restriction of $\md$ to $X_i$.
Since $X_i$ is connected and $\md_i$ is stable, $W_{\md_i}(X_i)$ is irreducible
of dimension $p_a(X_i)-1$, hence we are done (cf. \ref{Ssplit}).
\end{proof}

\begin{cor}
\label{irrth} $\Theta(X)$ is irreducible if and only 
if either $X$
is irreducible,
or every irreducible component of $X$ meets its complementary curve in
exactly $2$ points.
\end{cor}
\begin{proof} By \ref{thetadef}, $\Theta(X)$ is irreducible if and
only if $c=1$ (i.e. $X$ is free from separating nodes)  and
$\widetilde{b}=1$.

Assume $\Theta(X)$ irreducible; then $X$ has no separating nodes
and   $\widetilde{b} =\#\BX =1$. Hence  $\PXgb$ is irreducible, by \ref{Pstr}. 
Applying Corollary~\ref{irrP} we are done.

Conversely, if $X$ is irreducible, then $\Theta(X)$ is irreducible by
Theorem~\ref{irr}. 
If $X$ is reducible and satisfies the hypothesis,  obviously $c=1$.
Moreover, arguing as in the proof
of Corollary~\ref{irrP} we obtain that 
 $X$ has only one stable multidegree:
$\md= (\ov{g}_1,\ldots,
\ov{g}_{\gamma}) $, hence
$\Theta(X)$ is irreducible.
\end{proof}

\begin{remark}
In combinatorial terms, let
 ${\Gamma} _X^*$  be the graph obtained from $\Gamma_X$ by removing every
loop. Then
$\Theta(X)$ is irreducible if and only if either
 ${\Gamma}_X^*$ is a point, or
every vertex of ${\Gamma}_X^*$
has valency (i.e. degree)  $2$.
\end{remark}

\begin{remark}
\label{connect} Definition~\ref{thetadef}  coincides with the one given in
\cite{esttheta} or (which is the same) in
\cite{alex}, by    Theorem~\ref{Tstr}. In particular  $\Theta(X)$ is Cartier and ample.
\end{remark}

For the following simple Lemma  we use  the notation in
\ref{blow}.
\begin{lemma}
\label{key}
Let 
$S\subset \sing$ and   $M\in \Pic\YS$. 
Pick $\hM\in \Pic\hXS$ such that $\hM_{|\YS}=M$
and $\hM_{E}=\O_{E}(1)$ for every exceptional component $E$ of $\hXS$. Then
$
h^0(\hXS, \hM)=h^0(\YS,M).
$
\end{lemma}

\begin{proof}
(Cf. \cite{pacini} 2.1 for an analogous statement.)
For any pair of points $p_1,p_2\in \pr{1}$
choose a trivialization of $\O_{\pr{1}}(1)$ locally at such points;
now for
 any pair   
$a_1,a_2\in k$ there exists a unique section $s\in H^0( \pr{1},\O_{\pr{1}}(1))$
such that
$s(p_1)=a_i$ for $i=1,2$.
So, every section   $s_Y\in H^0(Y,M)$ extends to a unique section of $H^0(\hXS,\hM)$ determined by $s_Y$ and 
by the gluing data defining $\hM$.
Conversely, of course    any section in $H^0(\hXS,\hM)$  restricts to a  section of $M$.
\end{proof}

\begin{thm}
\label{Tstr} Let $X$ be a connected nodal curve. The stratification of $\PXgb$ given by \ref{Pstr}
induces the following canonical stratification:
\begin{equation}
\label{T1}
\Theta(X)=\coprod_{\stackrel{\emptyset\subseteq S \subseteq\sing}{\md \in \BYS}}\TSd,
\  \  \   \text{with  canonical isomorphisms}\  \  \TSd \cong W_{\md}(\YS).
\end{equation}
Equivalently, $
\Theta(X)=\{[M,S]\in \PXgb: \  h^0(\hXS,\hM)\neq 0\}.
$
\end{thm}
\begin{proof}
The equivalence of the two descriptions follows immediately from  \ref{Pstr} and Lemma~\ref{key}.
Furthermore it is clear that 
$$
\Theta(X)\subset \{[M,S]\in \PXgb: \ 
h^0(\hXS,\hM)\neq 0\}
$$
 (by uppersemicontinuity of $h^0$).
So we need to prove the other inclusion.

\

\noindent
{\bf Part 1. } { \it Proof  assuming $X$  free from separating nodes.}
In this case, by definition 
$$
\Theta(X)=\overline{\bigcup_{\md \in \Sigma({X})}\Wmd}.
$$
We shall use Abel maps (see \ref{abel}): recall that $\alpha^{\md}_{\YS}$ is the $\md$-th Abel map of
$\YS$ and the closure of its image in $\Pic^{\md}\YS$ is denoted by $A_{\md}(\YS)$.

\noindent 
{\it Step 1.} {\it Assume $\#S=1$ and let $\md \in \BYSs$. Then there exists $\me\in \BX$ such that
(using Theorem~\ref{irr} for the equality  below)
$$
\eSd(A_{\md}(\YS))\subset
\overline{\epsilon^{\me}_{\emptyset}(W_{\me}(X)})=
\overline{\epsilon^{\me}_{\emptyset}(A_{\me}(X))}.
$$
In particular, if $[M,S]\in \PXgb$ (so that $\mdeg M\in \BYS$)
satisfies
$\#S=1$ 
and $h^0(\hXS,\hM)\neq 0$, then $[M,S]\in \Theta(X)$.}

\noindent
Let $M\in \Pic^{\md}(\YS)$ with $M\in A_{\md}(\YS)$ and $\mdeg M=\md\in \BYSs$.
As $X$ is free from separating nodes, $\YS$ is connected.

Observe that, by \ref{Pstr}(\ref{smooth}), $\PXgb$ is the closure of its open subset
$$
\PXg= \coprod_{\me \in \Sigma({X})}P^{\me}_{\emptyset} \cong \coprod_{\me \in \Sigma({X})}\picX{\me}.
$$
Therefore there exists an $\me \in \BX$ such that $\eSd(M)\in \overline{P^{\me}_{\emptyset}}=
\overline{\picX{\me}}$.

Since $\#S=1$,
$|\md|=p_a(\YS)-1=g-2 = |\me|-1$.
Furthermore $\md \leq \me$
(by \ref{Pstr} (\ref{stcl})). Therefore there exists a unique index $i\in \{1,\ldots, \gamma\}$, say
$i=1$, such that $d_1=e_1-1$  and $d_i=e_i$ for $i\geq 2$.

Set  $S=\{n\}$, consider $\nu_S:\YS \to X$  the normalization at $n$ and let $C_1$ be the first
component of
$\YS$. Since $d_1=e_1-1$, by \ref{Pstr} (\ref{stcl})
$C_1$ contains one of the two branches of $n$,  call    $q_1$ this branch.
Let now $p_t\in C_1$ be a moving point specializing to $q_1$.

We can assume that $M$ is a general point in $ A_{\md}(\YS)$ (which is irreducible of  dimension
$p_a(\YS)-1$)  in particular that $M$ is in the image of the Abel map, that
$h^0(\YS,M)=1$, and that
$M$ has no base point lying over $n$ (by \ref{dimA} and \ref{bpf}).
Therefore there exists 
$L\in
\Pic X$ such that $\nu_S^*L=M$ and $L\in \im \amd$
(by \ref{h1} (\ref{h1gen})).
Set $L_t:=L(p_t)$;
then
$$
\mdeg L_t=\md + (1,0,\ldots,0)=\me \in \BX
$$
and 
$L_t\in \im \alpha^{\me}_X$, in particular
$h^0(X,L_t)\neq 0$. As $p_t$ specializes to $q_1$, we have that $\epsilon^{\me}_{\emptyset}(L_t)$
specializes to
$\eSd(M)$ so we are done with Step 1.

\

\noindent 
{\it Step 2.} {\it For every $S$ such that $\#S\geq 2$ and $\md \in \BYSs$, there exist 
$S'\subset S$ such that $\#S'=\#S -1$, and  a $\md ' \in \Sigma^{ss}(Y_{S'})$ such that}
$$
\eSd(A_{\md}(\YS))\subset\overline{\epsilon^{\md'}_{S'}(A_{\md '}(Y_{S'})}.
$$

Let $\md$ be a semistable multidegree for $\YS$. Consider the dual graph $\Gamma_{\YS}$
and an orientation on it inducing $\md$.
Note that $\Gamma_{\YS}$ is the subgraph of $\Gamma _X$ obtained by removing the edges corresponding to
$S$.
It is clear that if we add to $\Gamma_{\YS}$ any  edge  $n$ of $\Gamma$, 
(so that $n\in S$) oriented however we like,
we obtain a new oriented graph $\Gamma '$ such that 
$ \Gamma_{\YS}\subset \Gamma '\subset \Gamma_X.
$ 
Set $S'=S\smallsetminus \{n\}$, thus  $\Gamma '$ is the dual graph of the curve $Y_{S'}$
obtained by normalizing $X$ at $S'$.
Thus we have a map
$
\YS\to Y_{S'}
$
which is the normalization of $Y_{S'}$ at $n$.

The given orientation on $\Gamma '$ corresponds to a semistable multidegree $\md '$ such that $|\md
'|=|\md |+1$ and $\md '\geq \md $.

From now on we can argue  as for Step 1, with $Y_{S'}$ playing the role of $X$.
More precisely, if $\YS$ is connected, then the argument is exactly the same:
start from a general $M\in A_{\md}(\YS)$ and construct a family of line bundles
$L_t=L(p_t)\in A_{\md'}(Y_{S'})$ such that $p_t$ is a smooth point of
$Y_{S'}$ specializing to $n$, and $L\in A_{\md}(Y_{S'})$ such that $L$ pulls back to $M$.
Then
$\epsilon^{\md'}_{S'}(L_t)$  specializes to $\eSd(M)$. 

If $\YS$ is not connected,
then the general $M\in A_{\md}(\YS)$ has $h^0(M)\geq 2$, and it has no base point over $n$ 
(by \ref{bpf}).
We now apply \ref{d1} to obtain $L\in \im\alpha^{\md}_{Y_{S'}}$ which pulls back to $M$. The rest of the
argument is the same as before.

This conludes the proof of Step 2.

\

\noindent 
{\it Step 3.} {\it End of the proof of Part 1.}
To prove the theorem, we  pick $[M,S]\in \PXgb$ such that $M\in W_{\md}(\YS)$;
since $\md$ is stable, we have that $W_{\md}(\YS)=A_{\md}(\YS)$ by \ref{irr} (applied to every connected
component of $\YS$). 

Using Step 2  we can decrease
the cardinality of $S$ at the cost of passing from a stable multidegree to a semistable one
(which is why the assumption for Step 1 is that $\md$ is semistable, rather than stable).
Iterating Step 2 finitely many times,
we reduce the proof 
of the theorem to Step 1.
So the theorem is proved for $X$ free from separating nodes.

\

\noindent
{\bf Part 2.} {\it Proof assuming $\sep$ not empty.}
Recall that $\widetilde{X}\to X$ is the normalization of $X$ at $\sep$ and 
$\widetilde{X}=\cup_{i=1}^c X_i
$
denotes the decomposition of $\widetilde{X}$ into connected components; set $g_i=p_a(X_i)$.
By fact~\ref{Pstr} we have a canonical isomorphism
\begin{equation}
\label{Psplit}
\PXgb \cong \prod_{i=1}^c\overline{P^{g_i-1}_{X_i}}
\end{equation}
and, by Definition~\ref{thetadef}, another canonical isomorphism
\begin{equation}
\label{Tsplit}
\Theta(X) \cong \bigcup_{j=1}^c\bigr(\Theta (X_j) \times \prod_{\stackrel{i\neq
j}{1\leq i\leq c}}\overline{P^{g_i-1}_{X_i}}\bigl).
\end{equation}
Let $[M,S]\in \PXgb$ be such that $h^0(\YS, M)\neq 0$. Now $S\supset \sep$ hence we can factor
$$
\nS:\YS \stackrel{\widetilde{\nS}}{\la }\widetilde{X}\la X
$$
and denote $Y_i=\widetilde{\nS}^{-1}(X_i)$, so that $\YS$ is the disjoint union of $Y_1,\ldots Y_c$.
Note that $Y_i$ is the normalization of $X_i$ at a certain set of nodes, $S_i$,
of  $X_i$.
Therefore, under the isomorphism (\ref{Psplit}), the point $[M,S]$ corresponds to the point
$([M_1,S_1],\ldots, [M_c,S_c])\in \prod_{i=1}^c\overline{P^{g_i-1}_{X_i}}$
where $M_i=M_{Y_i}$.

Furthermore, $h^0(\YS, M)=\sum_1^ch^0(Y_i, M_i)$, hence there exists an index, say $i=1$, such that
$h^0(Y_1, M_1)\neq 0$.
Now, $X_1$ is free from separating nodes, therefore by the first part of the proof we obtain 
$ 
[M_1, S_1]\in \Theta(X_1).
$ 
By (\ref{Tsplit}), this implies $[M,S]\in \Theta(X)$ finishing the proof.
\end{proof}
\begin{example}
\label{ct} 
Let $X=C_1\cup C_2$ with $\#C_1\cap C_2=1$;
then $\BX$ is empty, while
$\Sigma (Y)=\{(g_1-1,g_2-1)\}$ ($Y$ is the normalization of $X$).
The points of $\PXgb$ correspond to line bundles 
of multidegree $(g_1-1,g_2-1)$ on $Y$ or to equivalence classes of line bundles
on
the curve $\hX$ obtained by blowing the unique node of $X$.
More precisely, if we order the components of $\hX$ so that  $\hX=C_1\cup E\cup C_2$
(where $E\cong \pr{1}$), then $\PXgb$ bijectively parametrizes line bundles
of multidegree $(g_1-1,1,g_2-1)$ on $\hX$.
There is a canonical isomorphism 
$$\PXgb\cong \Pic^{g_1-1}C_1\times \Pic^{g_2-1}C_2.$$

Now,  $\Theta(X)$ is canonically isomorphic to $W_{(g_1-1,g_2-1)}(Y)$, which we can easily
describe by means of \ref{sconn}.  
We obtain three different cases.

Case 1: $g_i\neq 0$ $i=1,2$. Then $\Theta(X)$ has two irreducible components:
\begin{equation}
\label{thetact}
\Theta(X)=\bigr( W_{g_1-1}(C_1)\times \Pic^{g_2-1}C_2\bigl)\cup
\bigr(\Pic^{g_1-1}C_1\times  W_{g_2-1}(C_2)\bigl).
\end{equation}
\end{example}
Case 2: $g_1=0$ and $g_2\neq 0$. Then  the first component  in (\ref{thetact}) is empty
and we get $\Theta(X)\cong W_{g_2-1}(C_2)\cong \Theta(C_2)$.

Case 3: $g_1=g_2=0$. Then $\Theta(X)$ is   empty.

\begin{example}
\label{theta2} Let $X=C_1\cup C_2$ with $\#C_1\cap C_2=\dd\geq 2$;  assume
$C_i$ smooth (this assumption can easily be removed) of genus $g_i$.  Then
$g-1=g_1+g_2+\dd -2$. We have
$ 
\BX=\{(g_1,g_2+\dd -2) ,(g_1+1,g_2+\dd -1),\ldots, (g_1+\dd -2,g_2) \},
$ 
so that $\PXgb$ has $\dd -1$ irreducible components of dimension $g$.
There is a canonical
isomorphism (cf.
\ref{Pstr} (\ref{smooth}))
$$
\PXg=\coprod_{i=0}^{\dd-2} P_{\emptyset}^{(g_1+i,g_2+\dd -i-2)}\cong \coprod_{i=0}^{\dd-2}
\Pic^{(g_1+i,g_2+\dd -i-2)}X.
$$
For every set $S\subset \sing$ such that $\#S=k$ with $1\leq k\leq \dd-2$,
we have 
$$
\Sigma (\YS)=\{(g_1,g_2+\dd -k-2) ,\ldots, (g_1+\dd -k-2,g_2)\};
$$
so that $\PXgb$ has a total of $(\dd-k-1){\dd\choose k}$ strata of codimension $k$,
each of which is isomorphic to $\Pic^{\md}\YS$.
If $k=\dd -1$ then for any choice of $\dd -1$ nodes, the curve obtained by blowing up $X$ at such nodes
has a separating node, hence $\BYS$ is empty.
Finally the last stratum corresponds to $S=\sing$ and $\md =(g_1-1,g_2-1)$, and it has codimension $\dd-1$.
We have
$$
P^{(g_1-1,g_2-1)}_{\sing}\cong \Pic^{g_1-1}C_1\times \Pic^{g_2-1}C_2.
$$

Now,  $\Theta(X)$ contains $\dd-1$  irreducible strata of dimension $g-1$,
one for every component of $\PXgb$.
Indeed for every $\md \in \BX$ we have
$\Theta^{\md}_{\emptyset}\cong \Wmd$ which is irreducible of dimension $g-1$,
by Theorem~\ref{irr}.

For every set $S\subset \sing$ such that $\#S=k$ with $1\leq k\leq \dd-2$,
  $\YS$ is connected and free from separating nodes,
so that for every $\md \in \BYS$ we get an irreducible stratum of dimension $g-k-1$
isomorphic to $W_{\md}(\YS)$.
If $k=\dd -1$ there are no strata (as before).
If $k=\dd$ we get a stratum isomorphic to the theta divisor
computed in Example~\ref{ct} (cf. (\ref{thetact})).
\end{example}

\section{Characterizing hyperelliptic stable curves}
We conclude the paper with a characterization of hyperelliptic irreducible curves, 
Theorem~\ref{W1thm}, extending a well known one for smooth
curves. The irreduciblity assumption is truly needed, as shown in counterexample~\ref{nonhyp}
\subsection{Irreducible  curves}
If we restrict  our interest to  irreducible singular curves, not only does the description of the compactified
jacobian simplifies substantially, but also,    the same description is valid for all degrees.

\begin{nota}{}
\label{descirr}
Let $X$ be an irreducible curve.  
Then the definitions of stable and semistable multidegrees (given for   $d=g-1$) coincide and are trivial. Thus,   for
every normalization $\YS\to X$ at a set $S$ of $\dS$ nodes, we have  $\BYS=\BYSs =\{p_a(\YS)-1\}=\{g-1-\dS\}$.
So, that definition   generalizes to all $d$, as follows.
With the notation of \ref{blow},
 a  line bundle $\hM\in \Pic^d\hXS$ is stable if (1) and (2) hold:
(1) $\deg_{\YS} \hM = d-\dS$.
(2) $\deg_{E_i} \hM = 1$ for all  $i=1,\ldots ,\dS$.

The equivalence relation is
the same as for $d=g-1$: two stable line bundles on $\hXS$ are   equivalent iff their pull backs to $\YS$
coincide. 
An equivalence class   is thus uniquely determined by
$S$ and by the restriction, $M$, of $\hM$ to $\YS$; we shall mantain the notation of \ref{blow}
and \ref{notpt}.

Exactly as in the case $d=g-1$,  we have the following.
{\it The variety $\PXb$ is reduced and irreducible. It
bijectively parametrizes  equivalence classes of stable line bundles on  the curves
$\hXS$  associated to $X$ as $S$ varies among all subsets of $\sing$.}
 
 Moreover, as in \ref{Pstr},
$\PXb$ has a canonical stratification in disjoint strata, called $P_S$,
indexed by the subsets $S$ of $\sing$. Every $P_S$ has a
canonical isomorphism (usually viewed as an identification)
$
\epsilon_S: \Pic^{d-\dS}Y_S\stackrel{\cong}{\la}P_S\subset \PXb.
$
 We have
\begin{equation}
\label{strirr}
\PXb=\coprod_{S\subset \sing}P_S\cong\coprod_{S\subset \sing}\Pic^{d-\dS}Y_S.
\end{equation}

%The stratification (\ref{strirr}) of $\PXb$ induces the  canonical stratification
%$$\Wdb\cong\coprod_{S\subset \sing}W_{d-\dS}(Y_S)$$
\end{nota}
\begin{nota}
\label{fam}
Given a family of  irreducible curves, $f:\X\to B$,
up to a finite  base change there exists the compactified Picard scheme $\pi_d:\pfb\to B$
which contains the relative degree-$d$ Picard scheme of $f$,
denoted $\Pic_f^d$,
as an open subset (see \cite{cner} for details).  The fiber of $\pi_d$ over a point $b\in B$ is $\ov{P^d_{X_b}}$.
\end{nota}
In the next Lemma we use the notation of \ref{notW}, in particular  (\ref{WMr}).
 \begin{lemma}
\label{W1} Let $\nu:\YS\to X$ be the   normalization of $X$ at a nonseparating  node $n$ of $X$,
set $\nu ^{-1}(n)=\{q_1,q_2\}$.
Let $M\in
W_{d}^r(\YS)$;
% so that $h^0(Y,M)\geq r+1$. 
then
\begin{enumerate}[(1)]
\item
\label{W1-}
$\WMr = \emptyset$ iff $h^0(\YS,M)=r+1$ and one of the two cases below occur:
 \begin{enumerate}[(a)]
\item
\label{a}
either   $ h^0(\YS,M-q_1-q_2))= h^0(\YS,M)-2,$ 
\item
\label{b}
or, up to
interchanging
$q_1$ with
$q_2$, 
$$
h^0(\YS,M)=h^0(\YS,M-q_1)\neq h^0(\YS,M-q_2).
$$
\end{enumerate}
\item
\label{W10}
$\dim \WMr = 0$ iff  $h^0(\YS,M)=r+1$ and
$$     h^0(\YS,M-q_1-q_2)=h^0(\YS,M -q_h)=r,\  \  h=1,2 .$$
In this case  $\WMr =\{L_M\}$ with $h^0(X,L_M)=r+1$.
\item
\label{W11}
$\dim \WMr = 1 $ iff one of the two cases below occur:
 \begin{enumerate}[(a)]
\item
\label{A}
$h^0(\YS,M)=h^0(\YS,M(-q_h))$ for $  h=1,2$ 
\item
\label{B}
$h^0(Y,M)\geq r+2$,
\end{enumerate}
\end{enumerate}
\end{lemma}
\begin{proof} It is a straightforward consequence of Lemma~\ref{d1}.
\end{proof}

\begin{nota}
\label{AC}
We recall a construction due to E.Arbarello and M.Cornalba  (cf. \cite{AC} section 2).
Let $h:T\to U$ be family of connected smooth  projective curves and 
assume  that $h$ has a section.
Then  for every pair of integers $(d,r)$,
there exists a $U$-scheme $\rho: W^r_{d,h}\to U$ 
such that for every $u\in U$, the fiber of $\rho$ over $u$ is
the Brill-Noether variety $W^r_d(h^{-1}(u))$ of the corresponding fiber of $h$.
Moreover there is a natural injective morphism of $U$-schemes,
$W^r_{d,h}\ha \Pic^d_h$,   viewed here as an inclusion.

Now let $f:\X\to B$ be a one-parameter family of 
smooth curves specializing to an irreducible curve  $X$,  let $b_0\in B$ be the
point over which the fiber   is $X$, and assume that the restriction of $f$ to $U=B\smallsetminus b_0$ is
smooth. Up to making a finite \'etale base change, we may asume that $f$ has a section 
(this  will not affect our conclusion). Call $h$ the restriction of $f$ to $U$ and introduce the
scheme
$W^r_{d,h}\to U$ described above. Consider the Picard scheme $\Pic^d_f\to B$, which is  integral,
separated and of finite type. Denote 
$\ov{W^r_{d,h}}\subset \Pic^d_f$ the closure of $W^r_{d,h}$ in $\Pic^d_f$.
Thus $\ov{W^r_{d,h}}$ is a scheme over $B$; we call 
$ W^r_{d,X}:=\ov{W^r_{d,h}}\cap \Pic^dX 
$ 
the fiber over $b_0$.
Then, by uppersemicontinuity of $h^0$, we have $W^r_{d,X}\subset W^r_d(X)$.
Therefore, 
{\it if $X$ is the specialization of a family of smooth curves $X_b$ such that
every irreducible component of $W^r_d(X_b)$ has dimension at least $c$ (for some number $c$),
%then every irreducible component of $ W^r_{d,X}$ has dimension at least $c$, and hence
then $\dim W^r_d(X) \geq c$} (i.e. $W^r_d(X)$ has a component of dimension at least $c$).
In particular:
{\it If $r\geq d-g$,
then  
$\dim W^r_d(X) \geq \BN=g-(r+1)(r-d+g)$.}
\end{nota}

\subsection{Hyperelliptic stable curves}
Some of the subsequent results are probably known to the experts, but an exhaustive reference was not found.

Let $H_g\subset M_g$  be  the locus of smooth
hyperelliptic curves   and $\ov{H_g}$ its closure in $\mgbar$.
We call a singular curve $X$      {\it hyperelliptic} if  it is
contained in
$\ov{H_g}$ (cf. \cite{HM}). 
%If $X$ is irreducible, \ref{hi} below applies, so there exists a unique
%line bundle $H_X$ in $W^1_2(X)$; we  call $H_X$ it the {\it hyperelliptic class} of $X$.

Some parts of the following proposition  can be found in, or easily derived from, \cite{CH} and \cite {HM}.
We here need a unified statement.
\begin{prop}
\label{hi}
Let $X$ be an irreducible nodal curve of genus   $g\geq 3$ with $\delta $ nodes and $\nu:Y\to
X$
its normalization. For every node $n_j$  set $\nu^{-1}(n_j)=\{q_1^j,q_2^j\}$.
The  following are equivalent.
\begin{enumerate}[(i)]
\item
\label{hi1}
There exists a line bundle $H_X \in \Pic^2X$ such that $h^0(X,H_X)=2$.
\item
\label{himg}
$[X]\in \ov{H_g}\subset \mgbar$ (i.e. $X$ is hyperelliptic).
\item
\label{hifam} There exists a family of smooth hyperelliptic curves $X_t$
specializing to $X$ and such that the hyperelliptic class of $X_t$ specializes to a line bundle, $H_X$, on $X$.

\item
\label{hi3}
There exists a $g^1_2$,\   $\Lambda$, on $Y$ such that $q_1^j+q_2^j$
is a divisor in $\Lambda$ 
for every
 $j=1\ldots\dd$ (in particular, $h^0(Y,q_1^j+q_2^j)\geq 2$).

\end{enumerate}
If the above hold, for every $j=1\ldots\dd$ we have $\nu^*H_X=\O_Y(q_1^j+q_2^j)$  and
$\Lambda \subset {\mathbb{P}} (H^0(Y,q_1^j+q_2^j)^*)$. Furthermore 
$W^1_2(X)=\{H_X\}$;  $H_X$ will be called  the {\emph {hyperelliptic class}} of $X$.
\end{prop}

\begin{remark}
\label{hired}
The implications  (\ref{hifam}) $\Leftrightarrow$ (\ref{himg}) and 
(\ref{hifam})$\Rightarrow$ (\ref{hi1}) hold also if  $X$ is reducible.\end{remark}

\begin{proof}
The implications (\ref{hifam})$\Rightarrow$ (\ref{hi1})  
and (\ref{hifam})$\Rightarrow$ (\ref{himg}) are obvious.  

 (\ref{hi1})$\Rightarrow $(\ref{hi3})
Let $\nu_1:Y_1\la X$ be the normalization of exactly one node
$n_1$ of
$X$. Let $M=\nu^*H_X$, then ($g_Y\geq 2$)  $h^0(Y_1,M)=2=h^0(X,H_X)$. Furthermore
$M$ is base-point-free, hence 
we are in
case (\ref{W10}) of Lemma~\ref{W1}.
We  obtain
  $M=\O_{Y_1}(q_1^1+q_2^1)$ and $H_X$ is uniquely determined
 (with the  notation of \ref{W1} (\ref{W10}),\  $H_X=L_M$).
Finally, $\Lambda_1 :={\mathbb{P}} (H^0(Y_1,M)^*)$; set $H_1=M$.

If $Y_1$ is smooth we are done, otherwise we iterate this procedure
as follows. 
Let $\nu_2:Y_2\to Y_1$ be the normalization of one node, $n_2$, of $Y_1$.
Call 
$$\nu_{1,2}:Y_2\stackrel{\nu_2}{\la} Y_1\stackrel{\nu_1}{\la} X
$$
and abuse the notation by using the same symbols for points in $X$, $Y_1$ and $Y_2$ whenever the normalization
maps are local isomorphisms.
 Then  $\nu_{1,2}^*H_X=\nu_2^*H_1=\nu_2^*\O_{Y_1}(q_1^1+q_2^1)=\O_{Y_2}(q_1^1+q_2^1)$.
Set $H_2=\nu_2^*H_1=\O_{Y_2}(q_1^1+q_2^1)$.
Note that the pull back
of the linear series $\Lambda _1$ to $Y_2$ is a $g^1_2$ containing $q_1^1+q_2^1$;
call $\Lambda_2=\nu_2^*\Lambda _1$ this $g^1_2$.
Now we distinguish two cases.

Case 1: $\dd \leq g-1$, i.e. $Y\neq \pr{1}$.

In this case we certainly have $p_a(Y_2)\geq 1$ 
hence $h^0(Y_2, H_2)=2$; thus we can argue as in the previous part to obtain $H_2=\O_{Y_2}(q_1^2+q_2^2)$
and $\Lambda_2 ={\mathbb{P}} (H^0(Y_2,q_1^j+q_2^j)^*)$ for $j=1,2$. This procedure can be repeated so we are
done.

Case 2: $Y= \pr{1}$.
We can argue as for case 1 only for $\dd -1$ steps,
  arriving at
$$
\nu:Y=\pr{1}\stackrel{\nu_{\dd}}{\la} Y_{\dd -1} \la X
$$
where $Y_{\dd -1}$ has only one node and all the above  morphisms are birational.
Furthermore,
for every $j=1,\ldots,\dd-1$ the pull back to $Y_{\dd -1}$ of $H_X$ is $\O_{Y_{\dd -1}}(q_1^j+q_2^j)$
 and $\Lambda_{\dd-1}={\mathbb{P}} (H^0(Y_{\dd-1},q_1^j+q_2^j)^*)$.

Now let $\Lambda:=\nu_{\dd}^*\Lambda_{\dd-1}\subset {\mathbb{P}} (H^0(Y,\O(2))^*)$.
For every $j=1,\ldots,\dd-1$ the divisor $q_1^j+q_2^j$ belongs to $\Lambda$ by construction.
To prove that also $q_1^\dd+q_2^\dd$ belongs to $\Lambda$ we repeat the same construction 
with respect to a different ordering of the nodes of $X$, for example by switching $n_{\dd}$ with $n_1$.
As $\Lambda $ is uniquely determined    by $H_X$,
and as $\dd \geq 3$, we are done.

(\ref{hi3})$\Rightarrow $(\ref{hi1})
Set $M=\O_Y(q_1^j+q_2^j)$ (for all $j$).
If $Y\neq \pr{1}$ we have $h^0(Y,M)=2$ 
and $h^0(Y,M-q_1^j-q_2^j)=1$,
so the proof is a straightforward
 iterated application of Lemma~\ref{W1}(\ref{W10}).

If $Y=\pr{1}$ we have $h^0(Y,M)=3$
and $M$ has no base point. Let $\nu_1:Y\to X_1$ be the map that glues only one pair of branches,
say $q_1^\dd,q_2^\dd$,
so that $p_a(X_1)=1$. Then for any $M_1\in \Pic X_1$ such that $\nu_1^*M_1=M$ we have
$h^0(X_1,M_1)=2$. 
Pick $M_1= \O_{X_1}(q_1^1+q_2^1)$ 
(abusing notation); we claim that for every $j=2,\ldots,\dd-1$ we have
$\O_{X_1}(q_1^j+q_2^j)\cong M_1$. 
This follows from the fact that, on $Y$,
the divisors $q_1^j+q_2^j$ belong all to the same $g^1_2$, $\Lambda$.
Indeed, recall that a line bundle on $X_1$ is uniquely determined by its pull back to $Y$, $M$, and by
the constant $c\in K^*$ gluing the two fibers  $M_{q_1^\dd}\stackrel{\cdot c}{\la}M_{q_2^\dd}$
via the multiplication by $c$. Furthermore, if $s\in H^0(Y,M)$ does not vanish at
$q_1^\dd$ and $q_2^\dd$, 
set $c(s)=s(q_2^\dd)/s(q_1^\dd)$, 
then $c(s)$ determines a unique line bundle $L_s $ which pulls back to $M$ and such that
the section $s$ descends to a section $\overline{s}\in H^0(X, L_s)$.
Now,  for every $j=1,\ldots \dd$, let $s_j\in H^0(Y,M)$ be   such that $div (s_j)= q_1^j+q_2^j$.
Then $M_1$ is uniquely determined by $c(s_1)$.
By hypothesis, the $\delta$ sections $s_j$ span a two dimensional subspace of $H^0(Y,M)$ and 
$s_\dd(q_1^\dd)=s_\dd(q_2^\dd)=0$; therefore we have that $c(s_j)=c(s_1)$ for every $j\leq \dd-1$,
proving that $\O_{X_1}(q_1^j+q_2^j)\cong M_1$ if $j\leq \dd -1$ (indeed, $div (\ov{s_j})=q_1^j+q_2^j$).

The claim enables us to complete the argument, again by lemma~\ref{W1} (\ref{W10}).

(\ref{himg})$\Rightarrow$ (\ref{hifam}).
If $X\in \ov{H_g}$ there exists a family of hyperellitic curves specializing to $X$.
Up to a finite base change, 
we get a family $f:\X\la B$ where $B$ is some smooth curve with a marked point $b_0\in B$,
such that the fiber $X_b$, $b\neq b_0$, is smooth and hyperelliptic, and the fiber over $b_0$
is
$X$. Moreover, 
 we get a line bundle $\mathcal H$ on $\X \smallsetminus X$
whose restriction to $X_b$ is the hyperelliptic line bundle on $X_b$.
The data $(\X\to B, {\mathcal H})$ induce a canonical map $\mu$ from $B\smallsetminus b_0$ to 
$\Pic^2_f$ such that $\mu(b)\in \Pic^2X_b$ is the   hyperelliptic
class of $X_b$ for all $b\in B\smallsetminus b_0$. As $B$ is a smooth curve $\mu$  extends to a regular map
$\mu:B\to \ov{P^2_f}$ (see \ref{fam}). 

We claim that $\mu(b_0)\in \Pic^2X\subset  \ov{P^2_X}\subset \ov{P^2_f}$.
By contradiction, suppose    $\mu(b_0)$ is a boundary point  of $\ov{P^2_X}$.
Therefore $\mu(b_0)=[M,S]$ where $S\subset \sing$ with $\dS=\#S\geq 1$ and 
$M\in \Pic^{2-\dS}\YS$.
 Since $\deg M\leq 1$ we   have   $h^0(\YS,M)\leq 1$. By lemma~\ref{key} we  get
    $h^0(\hXS,\hM)\leq 1$ for any representative $\hM$ for $[M,S]$.
 But $\hM$ is the specialization of line bundles having $h^0\geq 2$,
so this is impossible. 
The claim is thus proved, and so is the implication (\ref{himg})$\Rightarrow$ (\ref{hifam}).

Finally, we prove that
(\ref{hi3}) $\Rightarrow$ (\ref{himg}).
Let us denote by $G\subset \mgbar$ the locus of curves  satisfying
(\ref{hi3}). We claim that $G$ is irreducible of dimension
$2g-\dd -1$. 
Assume first  $\dd \leq g-1$; then
 $G$ is the locus of irreducible curves $X$
with
$\dd$ nodes, such that on the normalization $Y$   we have
$
h^0(Y,q_1^j+q_2^j)=2
$
and  if $\dd =g-1$ we need to impose also $q_1^j+q_2^j\sim q_1^{j'}+q_2^{j'}$. 
Thus a curve in $G$ is determined by its
normalization $Y$
and by the choice of $\dd$ linearly equivalent divisors of degree $2$ on $Y$.
As $\dim H_{g-\dd}=2(g-\dd)-1$
we get
$ \dim G = \dim H_{g-\dd} +\dd= 2g-\dd -1. $
Moreover, $G$ is irreducible because so is $H_{g-\dd}$.
If $\dd =g$, i.e. $Y= \pr{1}$  an element in $G$ is determined
by a $g^1_2$ on $\pr{1}$ and 
by $\delta$ divisors in it, everything  up to automorphisms. This yields $\dim G = 2+\delta - 3=\dd -1$.

Now denote by $\Delta_0^{\dd}$ 
the closure in $\mgbar$ of the locus of irreducible curves with  $\dd$ nodes.
It is well known that 
$\codim_{\mgbar} \Delta_0^{\dd}=\dd$. Therefore 
$ 
\dim (\ov{H_g}\cap \Delta_0^{\dd})\geq 2g-1-\dd
$ 
(as $\dim \ov{H_g}=2g-1$). Note  that (\ref{himg})$\Rightarrow$(\ref{hi3})
(we proved
(\ref{himg})$\Rightarrow$(\ref{hifam})$\Rightarrow$(\ref{hi1})$\Rightarrow$(\ref{hi3})),
hence
$ 
\ov{H_g}\cap \Delta_0^{\dd}\subseteq G.
$ 
As $\dim \ov{H_g}\cap \Delta_0^{\dd}\geq \dim G$, this   inclusion is an equality
and we are done.
\end{proof}
The next lemma  is  easy to prove for smooth
curves (cf. \cite{ACGH} p.13); our proof of the generalization  is 
 elementary  and maybe known, we include it  for completeness. 
 \begin{lemma}
\label{hysp}
Let $X$ be a hyperelliptic irreducible curve of genus $g\geq 3$; let $d$ and  $r$
be such that
$2\leq d\leq g$ and $0< 2r\leq d$.
Then 
%%?$ W^r_d(X)$ is irreducible of dimension 
$\dim W^r_d(X)=d-2r$.
\end{lemma}
\begin{proof} By definition,
$X$ is the specialization of some family of smooth hyperelliptic curves. The variety $W^r_d(C)$ of a smooth hyperelliptic curve $C$ 
is irreducible   of dimension
$d-2r$.
Therefore, by
the  construction of \ref{AC}, we obtain that
  $W^r_d(X)$ has dimension at least $d-2r$. So, it suffices to prove that every component of
$W^r_d(X)$ has dimension at
most $d-2r$ and that there exists   one component for which equality holds.
Furthermore, using the ``residuation" isomorphism
\begin{equation}
\label{resid}
W^r_d(X)\stackrel{\cong}{\la} W^{g-d+r-1}_{2g-2-d}(X)\  ;\hskip.7in  L\mapsto K_X\otimes L^{-1}
\end{equation}
we can reduce ourselves to prove the result for $d\leq g-1$.

Consider the partial normalization $\nu_n:Y_n\to X$  of one node $n$ of $X$ and let
$
\rho_r:W^r_d(X)\to W^r_d(Y_n)
$ be the pull back map. By Proposition~\ref{hi}
we have $\nu^*H_X=H_{Y_n}=\O_{Y_n}(q_1+q_2)$, where $\nu_n^{-1}(n)=\{q_1,q_2\}$.

We use induction on $\delta$. Suppose  $\dd =1$, we omit the subscript $n$
(i.e. $Y=Y_n$); now $g_Y=g-1$ and 
 $Y$ is a smooth hyperelliptic curve. 
$W^r_d(Y)$ is irreducible of dimension $d-2r$.
Let $U\subset W^r_d(Y)$ be the open dense subset $U=W^r_d(Y)\smallsetminus W^{r+1}_d(Y)$.
Pick $M\in U$, 
then (\cite{ACGH} p.13)    $M=H_Y^{\otimes r}(\sum
_{i=1}^{d-2r}p_i)$ with $h^0(Y,p_i+p_j)=1$ for all $i\neq j$.
By  Lemma~\ref{W1}
$W^r_M(X)$ is either empty or a single point; more precisely,  
$W^r_M(X)$ is not empty   exactly when neither $q_1$ nor $q_2$ appear among the $p_i$
(as  $h^0(M-q_1-q_2)=h^0(M\otimes H_Y^{-1})=h^0(M)-1$).
In this case 
every $\nu(p_i)$ is a smooth point of $X$, which we call again $p_i$; observe that
$h^0(X,H_X^{\otimes r}(\sum
_{i=1}^{d-2r}p_i))=r+1$,
therefore we necessarily have
$W^r_M(X)=\{H_X^{\otimes r}(\sum
_{i=1}^{d-2r}p_i)\}$.
We conclude that $\rho_r$ dominates $U$; more precisely, $W^r_d(X)$ has a unique irreducible component  
of
dimension equal to $d-2r$ dominating $U$. We also obtained that $\rho_r^{-1}(U)$ consists
of line bundles of the form $H_X^r(\sum
_{i=1}^{d-2r}p_i)$ with
$h^0(X,p_i+p_j)=1$ for all $i\neq j$.

The complement $ W^{r+1}_d(Y)$ of $U$ has dimension $d-2r-2$ 
 and the generic fiber of $\rho_r$ over it is a
$k^*$. Hence $\dim \rho_r^{-1}( W^{r+1}_d(Y))=d-2r-1$,
  so we are done.

Now assume $\dd\geq 2$. By the induction hypothesis, $W^r_d(Y_n)$ is irreducible of dimension $d-2r$
and $ W^{r+1}_d(Y_n)$ is either empty or irreducible of dimension $d-2r-2$.
We proceed as for $\dd =1$; set 
$U=W^r_d(Y_n)\smallsetminus W^{r+1}_d(Y_n)$ so that $U$ is irreducible of dimension $d-2r$.
By what we proved before, $U$ contains a non empty open subset $U_n$  consisting
of line bundles $M$ of the form
$M=H_{Y_n}^{\otimes r}(\sum
_{i=1}^{d-2r}p_i)$ with $h^0(Y_n,p_i+p_j)=1$ for all $i\neq j$.
By a trivial dimension count
we can disregard $U\smallsetminus U_n$ and concentrate on $U_n$.

Let $U_n'\subset U_n$ be the open subset of $M$ having  neither $q_1$ nor $q_2$ as base points;
 by Lemma~\ref{W1},   $W_M^r(X)$ is a single point for every $M\in U_n'$,
and $W_M^r(X)=\emptyset$ if $M\not\in U_n'$.
Therefore $W_M^r(X)$ has a unique  irreducible component of dimension $d-2r$
dominating $U_n$.
The rest of the proof  is  the same as for $\dd -1$. 
\end{proof}
The next result is well known if $X$ is nonsingular. 
\begin{thm}
\label{W1thm}
Let $X$ be irreducible of genus $g\geq 3$. Then 
%%?$W_{g-1}^1(X)$ is of pure dimension and
\begin{displaymath}
\dim W_{g-1}^1(X)=\left\{ \begin{array}{l}
g-3 \   \text{ if }\  X  \text{ is hyperelliptic } \\
g-4 \   \text{ otherwise.} \\
\end{array}\right.
\end{displaymath}
\end{thm}
\begin{proof}
If $X$ is hyperelliptic  this is a special case of Lemma~\ref{hysp}, so we will assume   $X$   not
hyperelliptic.
Now, for every smooth curve $C$ of genus $g\geq 3$, every irreducible component of $W^1_{g-1}(C)$ has dimension at
least $ g-4$   (and equality holds if and only if $C$ is not hyperelliptic).
Therefore  by \ref{AC},
$\dim W_{g-1}^1(X)\geq g-4$, hence
 it suffices to prove that 
\begin{equation}
\label{goal}
\dim W_{g-1}^1(X)\leq g-4
\end{equation}
(i.e.   every irreducible component   has dimension at most $g-4$).

If $g=3$ we are claiming that $W^1_2(X)$ is empty; this follows immediately from Proposition~\ref{hi}
(namely, from the fact that if $W^1_2(X)\neq \emptyset$ then $X$ is hyperelliptic).
So we shall assume $g\geq 4$ from now on.
Since $X$ is not hyperelliptic, by Proposition~\ref{hi} there exists a node $n$ of $X$ such that, denoting by
$\nu:Y\to X$ the normalization of $X$ at only $n$ and 
$\{q_1,q_2\}=\nu^{-1}(n)$,  we have
\begin{equation}
\label{nhyp}
h^0(Y, q_1+q_2)=1.
\end{equation}
Let us fix such a normalization, denote by $g_Y=g-1$ the genus of $Y$ and consider
the pull-back map
$$
\rho_1:W_{g-1}^1(X)\la W_{g-1}^1(Y)=W_{g_Y}^1(Y)
$$
defined by $\rho_1(L)=\nu^*L$. 
Recall that  $W_{g_Y}^1(Y)\cong W_{g_Y-2}^0(Y)$ (by (\ref{resid}))
hence
\begin{equation}
\label{1g}
\dim W_{g_Y}^1(Y)=\dim W_{g_Y-2}^0(Y)=g_Y-2=g-3.
\end{equation}
The fibers of $\rho_1$ have obviously dimension at most $1$.
Set 
$ \im\rho_1=I_0\dot{\cup} I_1$
where 
$$I_j=\{M\in \im\rho_1:\  \dim\rho_1^{-1}(M)=j\},\  \  \  j=0,1.$$ 
We shall prove (\ref{goal}) by showing that
\begin{equation}
\label{goal0}
\dim I_0\leq g-4
\end{equation}
and
\begin{equation}
\label{goal1}
\dim I_1\leq g-5.
\end{equation}
To prove (\ref{goal0}) we begin by observing that (\ref{nhyp}) is equivalent to
\begin{equation}
\label{nhypSD}
h^0(Y, \omega_Y(-q_1-q_2))=g_Y-2.
\end{equation}
Now it is easy to check that there exists a dense open subset $U\subset W^0_{g_Y-2}(Y)$ such that $\forall N\in U$
we have
$h^0(Y,\omega_Y(-q_1-q_2)\otimes N^{-1})=0$
(using \ref{00}).  Equivalently 
\begin{equation}
\label{SD2}
h^0(Y, N(q_1+q_2))=1,\  \   \forall N\in U.
\end{equation}
This implies that the  map $u$ below
\begin{equation}
%\label{}
\begin{array}{lccr}
u: & W^0_{g_Y-2}(Y)&\la &  \Pic^{g_Y}Y\\
&N &\mapsto &N(q_1+q_2)  
\end{array}
\end{equation}
satisfies
\begin{equation}
\label{dimim}
\dim \Bigr( u(W^0_{g_Y-2}(Y))\cap W^1_{g_Y}(Y)\Bigl)< \dim W^0_{g_Y-2}(Y)=g_Y-2.
\end{equation}
Now 
by Lemma~\ref{W1} we have
\begin{equation}
\label{I0}
I_0=\{M\in W^1_{g_Y}(Y): h^0(M-q_1-q_2)=h^0(M-q_h)=1,\   \  h=1,2\}.
\end{equation}
Therefore  $I_0\subset u(W^0_{g_Y-2}(Y))\cap W^1_{g_Y}(Y)$;  by (\ref{dimim}) we obtain
$
\dim I_0\leq g_Y-3=g-4
$
proving (\ref{goal0}).
To prove (\ref{goal1}) we apply Lemma~\ref{W1} to get
\begin{equation}
\label{I1}
I_1=\{M\in W^1_{g_Y}(Y): h^0(M-q_1)=h^0(M-q_2)=h^0(M)\}\cup  W^2_{g_Y}(Y);
\end{equation}
so we set $I_1 =J_a\cup J_b$
where $J_a:=\{M: h^0(M-q_h)=h^0(M)\geq 2, \   h=1,2\}$
and $J_b:=W^2_{g_Y}(Y)$.

The residuation isomorphism (\ref{resid}) gives
\begin{equation}
\label{resid1}
W^1_{g_Y-2}(Y)\cong W^2_{g_Y}(Y) = J_b.
\end{equation}
Assume   $Y$ is hyperelliptic, then by Lemma~\ref{hysp} we get
$ 
\dim J_b=g_Y-4=g-5
$ 
as wanted. Furthermore we have an injective
map
\begin{equation}
%\label{}
\begin{array}{lccr}
 & J_a&\ha &W^1_{g_Y-2}(Y) \\
&M &\mapsto &M(-q_1-q_2)
\end{array}
\end{equation}
hence, again by Lemma~\ref{hysp} we derive 
$ 
\dim J_a\leq \dim W^1_{g_Y-2}(Y)=g_Y-4=g-5
$ 
finishing the proof when $Y$ is hyperelliptic.
To conclude, observe that if (\ref{goal1}) holds in the special case of $Y$ hyperelliptic, it necessarily holds in the
generic case when $Y$ is not hyperelliptic, so we are done.
\end{proof}

\begin{example}
\label{nonhyp} The irreducibility hypothesis on $X$ cannot be removed from Theorem~\ref{W1thm}.
To see that, let $X=C_1\cup C_2$ be the union of two smooth curves meeting in one node $n$ of $X$;
call $q_i\in C_i$ the point corresponding to the node of $X$.
Recall that $X$ is hyperelliptic if and only if  $h^0(C_i,2q_i)=2$ for   $i=1,2$
(cf. \cite{CH}).

For any such $X$,  a description of   $\PXgb$ and of 
its theta divisor
has been given in Example~\ref{ct}. We   identify 
$ 
\PXgb=\Pic^{(g_1-1,g_2-1)}C_1\dot{\cup}C_2=\Pic^{g_1-1}C_1\times \Pic^{g_2-1}C_2
$ 
and  $\Theta(X)=\bigr( W_{g_1-1}(C_1)\times \Pic^{g_2-1}C_2\bigl)\cup
\bigr(\Pic^{g_1-1}C_1\times  W_{g_2-1}(C_2)\bigl)$.
Thus  we naturally  define
$$
W_{(g-1)}^{1}(X) =W_{(g_1-1,g_2-1)}^{1}(C_1\dot{\cup}C_2)\subset \Theta(X).
%=\{L\in \Pic^{(g_1-1,g_2-1)}C_1\dot{\cup}C_2: h^0(X,L)\geq 2\}
$$
Let us pick $C_1$  hyperelliptic of genus $g_1\geq 3$ and   $C_2$ non hyperelliptic
of genus $g_2\geq 3$. Hence $X$ is not
hyperelliptic.
Now we claim   $W_{(g-1)}^{1}(X)$ has a component of dimension $g-3$.
Indeed, consider
$W_{g_1-1}^1(C_1)\times \Pic^{g_2-1}C_2$. Since $C_1$ is hyperelliptic, $\dim W_{g_1-1}^1(C_1)=g_1-3$,
hence 
$$
\dim (W_{g_1-1}^1(C_1)\times \Pic^{g_2-1}C_2)=g_1-3+g_2=g-3.
$$
On the other hand,  it is clear that
$ 
W_{g_1-1}^1(C_1)\times \Pic^{g_2-1}C_2\subset W_{(g-1)}^{1}(X)
$ 
(indeed for every $M\in W_{g_1-1}^1(C_1)\times \Pic^{g_2-1}C_2$ we have $h^0(C_1\dot{\cup}C_2, M)\geq 2$).
\end{example}

\end{document}